%% file: mainArxiv.tex
\documentclass[12pt]{article}

\input{frontmatter.tex}

\begin{document}

\begin{center} 
\large Discrete Differential Geometry for $C^{1,1}$ Hyperbolic Surfaces of Non-Constant Curvature
\end{center}

\begin{center}
Christian Parkinson\footnote{Michigan State University, Department of Mathematics and Department of Computational Mathematics, Science and Engineering, 619 Red Cedar Rd, East Lansing, MI, 48824 (chparkin@msu.edu)} \hspace{2cm} Shankar Venkataramani\footnote{University of Arizona, Department of Mathematics, 617 N. Santa Rita Ave, Tucson, AZ, 85721 (shankar@arizona.edu)}
\end{center}

\noindent {\bf Abstract.} We develop a discrete differential geometry for surfaces of non-constant negative curvature, which can be used to model various phenomena from the growth of flower petals to marine invertebrate swimming. Specifically, we derive and numerically integrate a version of the classical Lelieuvre formulas that apply to immersions of $C^{1,1}$ hyperbolic surfaces of non-constant curvature. In contrast to the constant curvature case, these formulas do not provide an explicit method for constructing an immersion but rather describe an immersion via an implicit set of equations. We propose an iterative method for resolving these equations. Because we are interested in scenarios where the curvature is a function of the intrinsic material coordinates, in particular, on the geodesic distance from an origin or from an edge, we suggest a fast marching method for computing geodesic distance on manifolds. We apply our methods to generate surfaces of non-constant curvature and demonstrate how one can introduce branch points to account for the multi-generational buckling and subwrinkling observed in many applications. \\

\section{Introduction}

Thin elastic sheets that bend, buckle, and wrinkle arise in a variety of physical phenomena, both natural and man-made. Some applications where these appear are biological processes like leafy growth or marine flatworm locomotion \cite{Plant2,Plant1, Yamamoto}, the compression of steel during automotive collisions \cite{Steel1,Steel3,Steel2}, the irregular crenellations at the edges of torn plastic \cite{Sharon2004Leaves,Plastic1}, and photo- and thermo-sensitive hydrogels \cite{Hydrogel1,Hydrogel2}. 

\begin{figure}[b!]
\centering
\includegraphics[width = 0.45\textwidth]{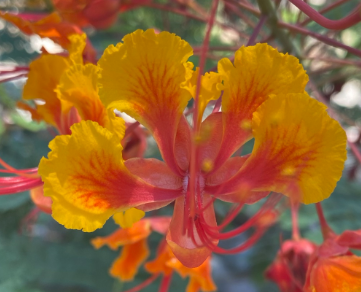} \, \includegraphics[width = 0.45\textwidth]{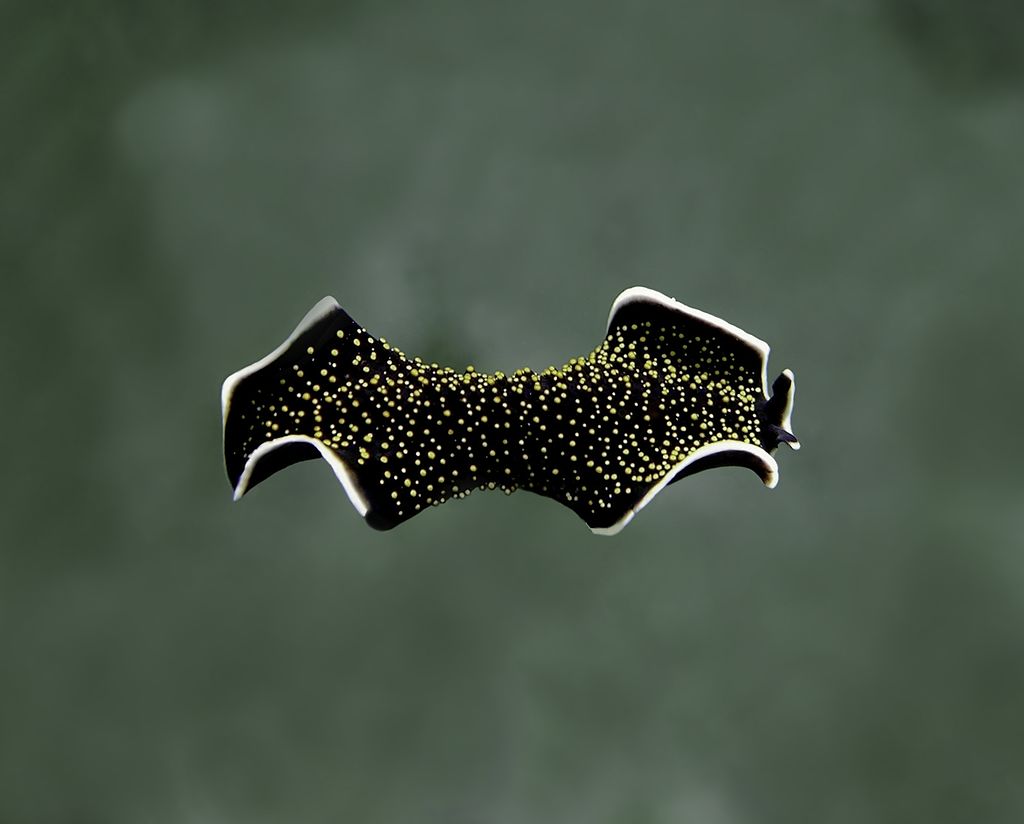}
\caption[caption]{Left: Red bird of paradise flower. The petals become curved and wrinkled near the edges. (Photo by CP)
 Right: Yellow papillae marine flatworm. While swimming, the center of the worm remains approximately flat while the edges become curved. (Photo by Betty Wills, uploaded to Wikipedia Commons under Share-Alike license: \url{https://commons.wikimedia.org/wiki/File:Yellow\_papillae\_flatworm\_(Thysanozoon\_nigropapillosum).jpg}) }
\label{fig:rbop}
\end{figure}

Over the last decade, the theory of non-Euclidean elastic plates has been developed to model such phenomena \cite{Elastic1,Elastic2}. A key underlying assumption is that the observed configurations of these plates are minimizers of an elastic energy functional \cite{Elastic1}. The goal is then to create immersions of bounded hyperbolic surfaces into $\R^3$ that minimize elastic energy, and the analysis of these immersions is central to some recent work devoted to generating these surfaces \cite{Intro3,Intro1,Intro4}. For an overview of the calculus of variations for thin elastic sheets, we direct the reader to \cite{Lewicka}.

In the limit of vanishing thickness, the sheet becomes essentially inextensible and the energy minimization problem enforces a hard isometry constraint on the immersion \cite{Elastic2}. 
In this manuscript, we will be concerned with modeling 2-dimensional Riemannian surfaces $(\Omega, g)$ of non-constant negative curvature realized as isometric immersions into $\R^3$. Mathematically, this amounts to solving a nonlinear PDE, the prescribed Gauss curvature equation \cite{stoker}, for an immersion $\mathbf{r}: \Omega \to \mathbb{R}^3$. This equation is a fully nonlinear hyperbolic Monge-Amp{\`e}re equation if the prescribed curvature is negative. 

A separate, yet related, line of inquiry in the field of geometric mechanics concerns the application of the Monge-Amp{\`e}re equation to prestrained plates with positive curvature, which are governed by the {\em elliptic} Monge-Amp{\`e}re equation. While the analytic theory for this equation is well-established, the development of corresponding computational methods is a more recent endeavor. Difficulties arise from several sources. Beyond the equation's inherent nonlinearity, weak solutions are not based on standard variational principles but on concepts such as viscosity solutions, which can be challenging to interpret discretely \cite{Nochetto2019TwoScale}. Furthermore, the global {\em convexity constraint}, often required for a well-posed theory, is difficult to enforce in a numerical setting. Consequently, the construction of convergent numerical schemes remains an active area of research, even for this well-studied, elliptic class of Monge-Amp{\`e}re equations \cite{NEILAN2020105}.

In contrast, even the analytic theories for {\em parabolic} and {\em hyperbolic Monge-Amp{\`e}re equations} are very much under active development. On the numerical side, there have been several recent developments in the simulation of the parabolic Monge-Amp{\`e}re equation \cite{BARTELS2020221}, which arises in applications such as {\em curved crease origami} in flat sheets \cite{Liu2024Design}. The numerical solution of this equation is complicated by the presence of curved creases, for which many standard methods fail to converge  \cite{bartels2025babuvskasparadox}. This highlights the distinct challenges posed by different classes of Monge-Amp{\`e}re equations and the motivation for developing specialized numerical techniques (See \cite{Bonito2020PartI, Bonito2021PartII} and references therein).

In this paper, we develop methods for the hyperbolic setting. 
In particular, motivated by the $C^1$ immersion theorem of Nash-Kuiper \cite{Kuiper,Nash} and the nonexistence of $C^2$ immersions for complete (unbounded) hyperbolic surfaces proven by Efimov \cite{Efimov,Milnor}, Shearman and Venkataramani \cite{Toby} developed the theory of distributed branch points for negatively curved surfaces, which allows immersions whose smoothness is strictly between $C^1$ and $C^2$. They argue that such immersions can have significantly lower elastic energy than $C^2$ immersions of finite pieces of hyperbolic surfaces, especially as these pieces become large, and can successfully capture the sort of multigenerational buckling and subwrinkling that appear in nature and are not captured by smoother immersions. They also develop discrete differential geometry (DDG) for such surfaces, and suggest a method for constructing discrete immersions by {\em surgery}, i.e. recursively patching together ``sectors'' --- bounded hyperbolic surfaces that admit a pair of {\em asymptotic coordinates} \cite{stoker} and have data specified along two transversely intersecting boundary curves 
\cite{Hartman1952Hyperbolic}.

It is worth emphasizing that the methods developed in \cite{Toby}, and extended in this manuscript, are intrinsically discrete rather than being discrete approximations of a continuous formulation. Our approach borrows extensively from the concepts of Discrete Differential Geometry (DDG), which aims to construct discrete theories that run parallel to the smooth theory while preserving its essential geometric and structural properties. This philosophy guides our development of the discrete Lelieuvre formulas in Section \ref{sec:discLeli}.

In general, one cannot expect to obtain solutions for a nonlinear PDE as a piecewise function obtained from local solutions. However, the prescribed Gauss curvature equation for immersions of hyperbolic surfaces has a special structure --- it is the compatibility condition \cite{BobenkoSuris} for the {\em Lelieuvre equations} \cite{eisenhart1997riemannian}, a pair of ODEs defining the immersion $\mathbf{r}$ along {\em asymptotic curves} \cite{hartman1951asymptotic-943,stoker} which are the characteristic curves for the underlying hyperbolic PDE. The key insight in \cite{Toby} was that if two solutions agree along an asymptotic curve, then they can be patched together along this curve to obtain a piecewise defined weak solution for an immersion with a prescribed Gauss curvature. This naturally leads to the question of obtaining solutions on sectors, i.e. regions bounded by a pair of asymptotic curves. Each sector is then obtained by solving a {\em Goursat problem}, i.e. a second order hyperbolic PDE with data specified along two characteristic curves \cite{Hartman1952Hyperbolic}, which, with appropriate choices for the boundary data, can be patched together with other sectors. In particular, this procedure can be carried out using sectors where one or both boundaries are straight lines in $\mathbb{R}^3$. Hyperbolic surfaces that contain two intersecting straight lines were first obtained by Amsler \cite{Amsler} and their discrete analogs have also been studied \cite{Bobenko2002Nonlinear}. Each of the four quadrants defined by these intersecting lines is an example of an {\em Amsler sector} \cite{Toby}. In addition, the surgery procedure also uses {\em Pseudo-Amsler} sectors -- sectors where one of the boundary curves is a straight line, but the other boundary is a prescribed space curve in $\mathbb{R}^3$ \cite{Toby}. 

All the work in \cite{Toby}---as well as that in \cite{Intro3,Intro1}---applies to surfaces of constant negative curvature. In this work, we extend these methods to surfaces of non-constant negative curvature. In physical scenarios, this could account for thin sheets that are relatively flat in some places but become increasingly curvy in others. We are specifically interested in scenarios where curvature depends on the geodesic distance from some ``center" or an ``edge''---which could accurately model objects like leaves or flower petals that become frilly only near the edges---though the methods described could be applied more generally. A key tool for us will be the Lelieuvre formulas \cite[Chap 6]{eisenhart1997riemannian}. These formulas establish a relationship between an immersion of the surface and the normal vector along the surface and provide a relatively simple way to generate hyperbolic surfaces with prescribed negative curvature by evolving the surface ``outward" from consistently defined boundary data. However, in the case of non-constant curvature, this introduces some difficulty in synthesizing the different coordinate systems used to describe the surface. The connection between these intrinsic and extrinsic coordinates is the immersion itself. This creates a dependency loop: the Lelieuvre formulas, which are expressed in extrinsic asymptotic coordinates, require the curvature as an input; however, the curvature is defined in terms of intrinsic material coordinates. Resolving this connection—constructing an immersion that depends on information that can only be fully known once the immersion itself is constructed—is a central challenge we overcome in this work.

The organization of this paper is as follows. In section \ref{sec:contLeli}, we present a self-contained derivation of the Lelieuvre formulas for a surface of non-constant negative curvature. In section \ref{sec:discLeli}, we develop the discrete analog to the continuous Lelieuvre formulas. In section \ref{sec:geodesic}, we describe a fast marching method for computing geodesic distance on triangulated manifolds, which will allow us to generate surfaces with prescribed negative curvature depending on the geodesic distance. In section \ref{sec:amsler}, we apply our methods to generate Amsler surfaces with non-constant negative curvature. In section \ref{sec:branch}, we describe how the method of ``surgery" developed by \cite{Toby} can be adapted to our scenario, and thus generate $C^{1,1}$ surfaces with branch points that allow for buckling and subwrinkling. While sections \ref{sec:amsler} and \ref{sec:branch} utilize a disk-like geometry, in section \ref{sec:strip_geometry}, we present a further application of our methods in a geometry which models the edge of a torn plastic strip as studied in \cite{Sharon2004Leaves,Plastic1}. Finally, in appendix \ref{sec:convProof}, we give a proof that our iterative algorithm for generating a discrete surface of prescribed, non-constant, negative Gauss curvature does indeed converge.

\section{Lelieuvre Formulas for Surfaces with Non-Constant Negative Curvature} \label{sec:contLeli}

In this section, we present a short and self-contained derivation of the Lelieuvre formulas \cite[Chap 6]{eisenhart1997riemannian}, a set of differential equations describing the geometry of negatively curved surfaces. Our thin elastic object is modeled by a 2-dimensional Riemannian manifold $(\Omega,g)$ representing its mid-surface. The metric tensor $g$ determines the intrinsic geometry on $\Omega$ \cite{stoker}, and is called the {\em target metric} in the context of non-Euclidean elasticity \cite{Elastic1}. 

Supposing our surface $\Omega$ has a pair of reference coordinates $(x,y)$, we begin with an immersion $\mathbf{r}: \Omega \to \mathbb R^3$ of our surface into Euclidean space. The immersion induces a metric $I$ on $\Omega$, called the {\em first fundamental form} \cite{stoker}, by pulling back the standard Euclidean metric on $\mathbb{R}^3$, i.e. $I = d\mathbf{r} \cdot d\mathbf{r}$. The difference $I-g$ between the induced and the target metrics on $\Omega$ is a measure of the stretching of the sheet. In the thin limit, sufficiently regular elastic objects are essentially unstretchable \cite{Elastic2}, so the allowed configurations satisfy the {\em isometry constraint} $d\mathbf{r}\cdot d\mathbf{r}=g$.

\begin{figure}[b!]
\centering
\includegraphics[width = 0.8\textwidth]{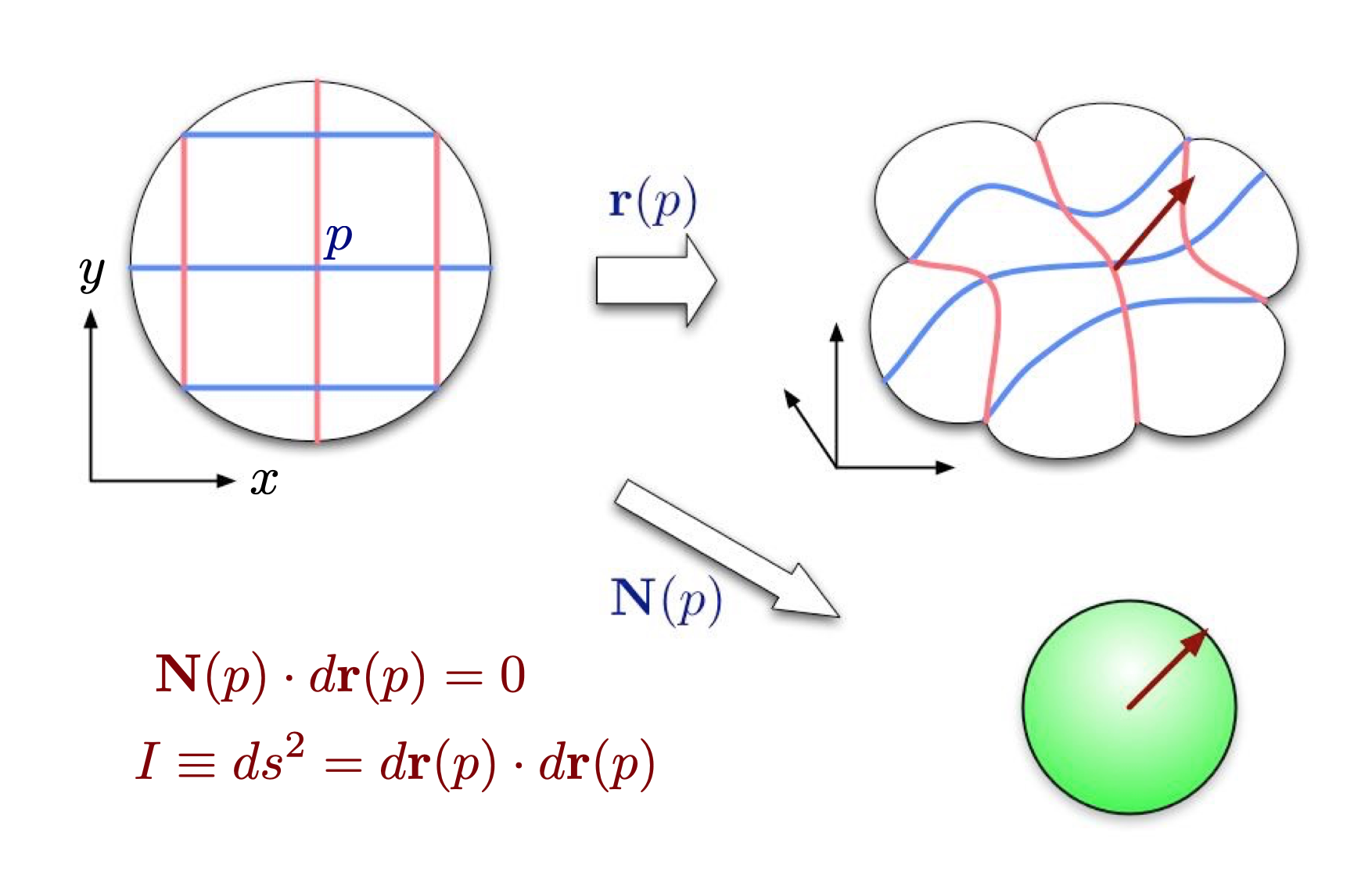}
\caption{A schematic representation of the mappings associated with the immersion $\mathbf{r}:\Omega \to \mathbb{R}^3$.}
\label{fig:geometry}
\end{figure}

Since $\br$ is an immersion, at any point $p$ on the surface, $\br_x$ and $\br_y$ form a basis for the tangent space $T_p$. In particular, $\br_x \times \br_y$ is never zero, and the unit normal to the surface can be written \begin{equation} \label{eq:normalDef}\mathbf{N} = \frac{\br_x \times \br_y}{\abs{\br_x \times \br_y}}. \end{equation} The map $\mathbf{N}: \Omega \to S^2$ given by $p \mapsto \mathbf{N}(p)$ is called the {\em Gauss normal map} \cite{stoker}.  Fig.~\ref{fig:geometry} shows a schematic representation of these mappings.

Note that for a $C^k$ surface  $\mathbf{r}(x,y)$ the normal map is only $C^{k-1}$.  If the normal map is differentiable at a point $p$, by differentiating $\abs{\bN} =1$, we see $\bN \cdot \bN_x = \bN \cdot \bN_y = 0$, which implies that $\bN_x$ and $\bN_y$ also lie in the tangent plane $T_p$. The \emph{shape operator} at the point $p$ is the unique linear operator $S_p : T_p \to T_p$ such that $S_{p}(\br_x(p)) = \bN_x(p)$ and $S_p(\br_y(p)) = \bN_y(p)$ \cite{stoker}.  The shape operator provides information regarding local curvature of the immersed surface $\br(\Omega)$ near $p$. The determinant of the shape operator defines the Gaussian curvature extrinsically: $K(p;\br):= \text{det}(S_p)$. From this, we see \begin{equation}\label{eq:TheoremaEgregium}
\bN_x \times \bN_y = \text{det}(S_p) \br_x \times \br_y = K(p;\br) \abs{\br_x \times \br_y} \bN
\end{equation} or $\abs{K(p;\br)} = \frac{\abs{\bN_x \times \bN_y}}{\abs{\br_x \times \br_y}}$ with the sign determined by the relative orientations of $\bN$ and $\br_x \times \br_y$. This states that the Gauss curvature is the local magnification of an infinitesimal area element on the surface under the normal map $\mathbf{N}: \Omega \to S^2$ given by $p \mapsto \mathbf{N}(p)$. We note that \emph{a priori} $K(p;\br)$ may also depend on the immersion $\br$. However, as a consequence of the Theorema Egregium, Gauss curvature is intrinsic \cite{stoker} which implies that if $d\mathbf{r} \cdot d \mathbf{r} = g$, then $K(p;\br)$ is determined by $g$ and independent of the particular immersion $\br$. This yields the {\em prescribed Gauss curvature} equation
\begin{equation}
    \text{det}(S_p) = K_g(p)
    \label{eq:Gauss}
\end{equation}
where the left hand side is the determinant of the shape operator, given in terms of the immersion $\mathbf{r}$ and the right hand side is given by a function $K_g:\Omega \to \mathbb{R}$, where $K_g$ is the curvature for the metric $g$. Equation ~\eqref{eq:Gauss} is a fully nonlinear second-order PDE for the immersion $\mathbf{r}$. In this work, we will numerically solve this equation to find solutions $\mathbf{r}:\Omega \to \mathbb{R}^3$ for hyperbolic surfaces, i.e., specifically for situations where $K_g(p) < 0$. We henceforth drop the subscript $g$. 

For a broader context, equation \eqref{eq:Gauss} is an example of a Monge-Amp\`ere equation, that is hyperbolic (resp. elliptic) if $K_g < 0$ (resp. $K_g > 0$). In contrast to the case of convex solutions of elliptic Monge-Ampere equations, that have a well-developed regularity theory \cite{Pogorelov1964Book,Caffarelli1990Localization}, the regularity theory for the hyperbolic case, and even the question of the appropriate definition of the notion of a weak solution, (when considered on sets of various different codimensions in spaces of various different dimensions) is very much an active area of research \cite{Cao2,Cao1,Cao3,LewickaMonge2,LewickaMonge3,LewickaMonge1}. In this work, we will be primarily concerned with immersions $\mathbf{r}:\Omega \to \R^3$ that are $C^{1,1}$ and piecewise smooth, but not necessarily globally $C^2$. Despite the lack of two derivatives for $\br$ it is straightforward to interpret such $C^{1,1}$ immersions as weak solutions of \eqref{eq:Gauss} by multiplying both sides of the equation by compactly supported smooth functions and integrating over the domain $\Omega$. Indeed, we only need a.e. defined measurable functions $\text{det}(S_p)$ and $K_g(p)$ for this procedure to make sense. In our setting $K_g$ is a ``fixed" smooth function, i.e. independent of the immersion $\mathbf{r}: \Omega \to \R^3$. For a $C^{1,1}$ immersion $\br$, the normal map $\mathbf{N}:\Omega \to S^2$ is Lipschitz, and Rademacher's theorem \cite[\S 3]{Heinonen2005Lipschitz} implies that $\mathbf{N}$ is a.e. differentiable, so the shape operator $S_p$ is defined a.e. 

If an immersion $\mathbf{r}:\Omega \to \R^3$ is smooth in a neighborhood of a point $p$ and the Gauss curvature is negative at a point $p$, then the shape operator $S_p$ is well defined at $p$ and has one negative and one positive eigenvalue, and thus there are two linearly independent vectors in $T_p$ each satisfying $\bv \cdot S_p(\bv) = 0$. These vectors define the {\em asymptotic directions} at $p$: the directions along which there is zero normal curvature \cite{stoker}. Assuming the Gaussian curvature is negative everywhere and the vectors defining these asymptotic directions vary smoothly with $p$, we can extend them to integral curves which give asymptotic coordinates on the manifold \cite{hartman1951asymptotic-943,stoker}. We denote the asymptotic coordinates by $(u,v)$. These coordinates are not intrinsic to the surface, i.e. tied to material points, but rather they depend on the particular immersion $\br: \Omega \to \mathbb{R}^3$. 

From the defining property of the asymptotic coordinates $(u,v)$  we get
\begin{equation}\label{eq:noNormalCurvature}
\begin{matrix}
\br_u \cdot S(\br_u) = 0, \\
\br_v \cdot S(\br_v) = 0,
\end{matrix}  \,\,\,\,\,\,\,\,\,\,\, \implies \,\,\,\,\,\,\,\,\,\,\, \begin{matrix}
\br_u \cdot \bN_u= 0, \\
\br_v \cdot \bN_v = 0.
\end{matrix}
\end{equation}

This implies that both $(\br_u, \bN, \bN_u)$ and $(\br_v, \bN, \bN_v)$ are orthogonal frames. In particular, \begin{equation}\label{eq:abdef}
\br_u = a\bN_u \times \bN, \,\,\,\,\,\,\,\,\,\, \br_v = b \bN_v \times \bN,
\end{equation}
for some functions $a(p)$, $b(p)$. We resolve $a(p)$ and $b(p)$ by equating the mixed partial derivatives obtained from either of the above expressions: \begin{equation}\label{eq:mixedPartialR}
\begin{aligned} \br_{uv} &= a_v\bN_u \times \bN + a \bN_{uv} \times \bN + a \bN_v \times \bN_{u}\\ &= b_u \bN_v\times\bN + b\bN_{uv} \times \bN + b \bN_{u}\times \bN_{v}. \end{aligned}
\end{equation}
The first and second terms in each expression in \eqref{eq:mixedPartialR} are normal to $\bN$, so for the equation to hold, the last terms must be equal which shows that $a(p) = -b(p)$. From \eqref{eq:TheoremaEgregium}, we see $a(p)b(p) = K(p)^{-1}$ so that \begin{equation}\label{eq:abK}
a(p) = -b(p) = (-K(p))^{-1/2}.
\end{equation} Exploiting $\bN \times \bN = 0$, we can write \begin{equation}\label{eq:lelieuvre1} \br_u = (-K(p)^{-1/2}) \bN_u \times \bN = [(-K(p)^{-1/4})\bN]_u \times [(-K(p)^{-1/4})\bN] \end{equation} and similarly in the equation for $\br_v$. Finally, defining the scaled normal $\bnu = (-K(p)^{-1/4})\bN$ \cite{schief2017gaussian-047}, we arrive at \begin{equation} \label{eq:lelieuvre} \br_u = \bnu_u \times \bnu, \,\,\,\,\,\,\,\,\,\, \br_v = -\bnu_v \times \bnu.\end{equation} Equating the mixed partial derivatives of $\br$ yields the compatibility condition \begin{equation}\label{eq:compatibility} \bnu \times \bnu_{uv} = 0. \end{equation} The formulas in \eqref{eq:lelieuvre} and \eqref{eq:compatibility} are a version of the Lelieuvre formulas \cite{eisenhart1997riemannian}  with non-constant curvature. For pseudospherical surfaces---those with constant negative curvature---one can choose units of length so that the curvature is $K(p) = -1$, which gives the classical Leliuvre formulas (and compatibility condition): \begin{equation} \label{eq:lelieuvreConstCurv} \begin{split} &\br_u = \bN_u \times \bN, \,\,\,\,\,\,\,\,\,\, \br_v = -\bN_v \times \bN, \\ &\,\,\,\,\,\,\,\,\,\,\,\,\,\,\,\,\,\,\,\,\,\bN \times \bN_{uv} = 0. \end{split} \end{equation} In this manuscript, we consider curvature functions which are nonconstant, but are negative and bounded away from zero. 

\section{Discrete Lelieuvre Formulas} \label{sec:discLeli}

\begin{figure}
\centering
    \includegraphics[width=0.45\textwidth,trim=40 20 35 20,clip]{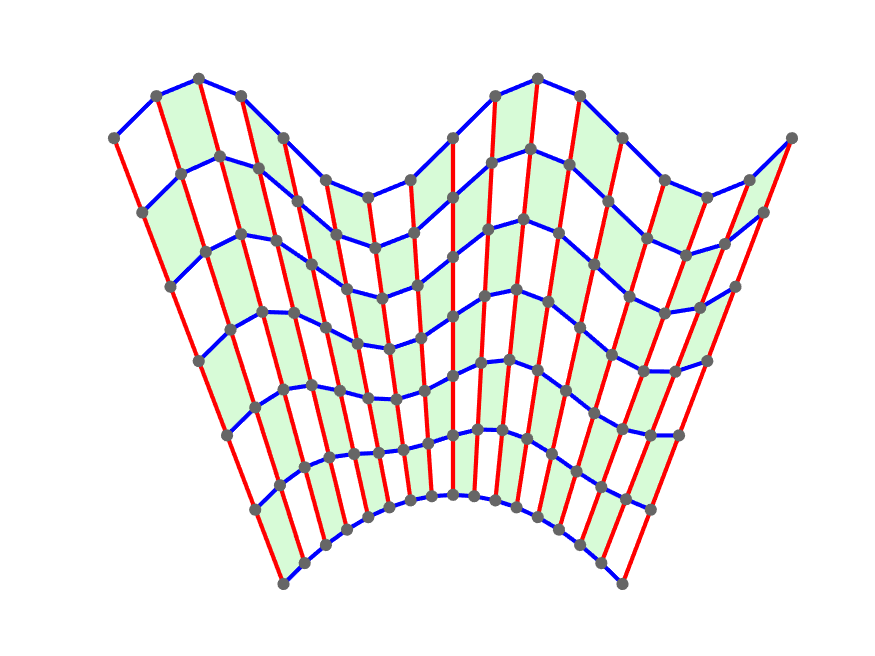} \,\,\,\,\,\,
    \includegraphics[width=0.45\textwidth,trim=40 20 35 20,clip]{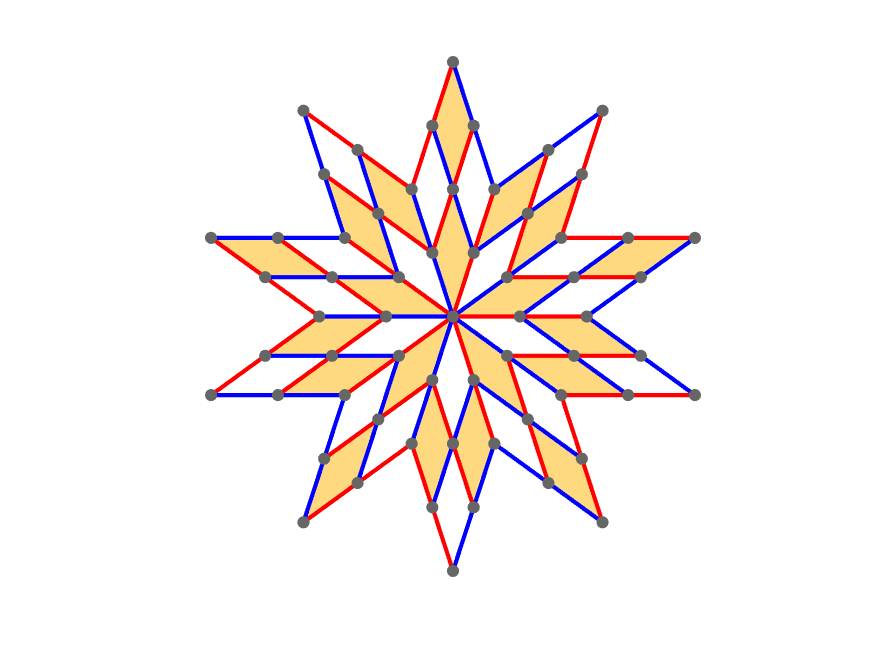}
    \caption{Two examples of asymptotic complexes projected onto the plane. For each, the degrees of all interior vertices are even, and this allows the complex to be ``checkered" \cite{Toby}. On the left, each interior vertex has degree 4, which allows the local asymptotic coordinates on each quad (represented by red and blue lines) to be extended globally. By contrast, the asymptotic complex on the right can be naturally interpreted as 10 sectors that have been pasted together as described in section \ref{sec:amsler}. Accordingly, the central vertex has degree 10, and asymptotic coordinates can be consistently extended to pairs of adjoining sectors, but not to the entire complex. This is discussed in section \ref{sec:amsler} and demonstrated in figure \ref{fig:patching}.}
    \label{fig:Acomps}
\end{figure}

In this section, we describe some of the underlying discrete differential geometry which we will use to construct discrete surfaces with prescribed negative curvature. The discrete object we use to approximate a negatively curved surface is an asymptotic complex as defined in \cite{Toby}. For our purposes the full definition is not crucial, and it suffices to say that a asymptotic complex is a strongly regular quadgraph (that is, a collection of vertices, edges and faces where every face is a quadrilateral and the intersection of two faces is either empty, a vertex, or an edge) where the set of edges can be decomposed into two disjoint subcollections such that edges bounding each face are taken from alternating subcollections. Two examples of asymptotic complexes are shown in figure \ref{fig:Acomps}. It is proven in \cite{Toby}, and demonstrated in figure \ref{fig:Acomps}, that all interior vertices in an asymptotic complex have even degree if and only if the asymptotic complex can be colored with two colors such that no faces which share an edge have the same color. Additional structure is added to the definition in \cite{Toby} that allows for some interfacing between the discrete complex and the continuous object that it approximates, which will not be necessary for our purposes. For background on cell decompositions of surfaces and related concepts, we refer the reader to \cite[Chap. 0]{Hatcher}, \cite{Huhnen}, \cite[\S2]{Kacz}.

Since $\R^2$ has the natural structure of an asymptotic complex with vertices given by $\mathbb Z^2$, it is simplest to focus on discrete immersions \begin{equation} \label{eq:discImmersion} \begin{split} &\br: \mathbb Z^2 \to \mathbb R^3,\\ &(i,j) \mapsto \br_{ij}, \end{split}\end{equation} that obey a discrete version of the Lelieuvre formulas. Here the quadrilaterals are the squares with coordinates $(i,j), (i+1,j), (i+1,j+1), (i,j+1)$ for some fixed $i,j \in \mathbb Z$, and we associate the discrete coordinate $i$ with the asymptotic coordinate $u$, and the discrete coordinate $j$ with the asymptotic coordinate $v$ (that is, the horizontal edge connecting $(i,j)$ to $(i+1,j)$ is the projection of a $u$-asymptotic line down to $\R^2$). 
We adopt a notation which is standard in discrete differential geometry \cite{BobenkoSuris}: for any quantity $f_{ij}$ defined on a quadrilateral, we write $$f_0 \defeq f_{ij}, \,\,\,\, f_1 \defeq f_{i+1,j}, \,\,\,\, f_2 \defeq f_{i,j+1}, \,\,\,\, f_{12} \defeq f_{i+1,j+1}.$$ Using this notation, we say an immersion $\br:\Z^2 \to \R^3$ is a discrete $K$-surface---that is, a discrete surface with prescribed Gaussian curvature $K$---if there is a discrete normal field $\bN:\Z^2 \to S^2$ such that \begin{equation} \label{eq:discreteLelieuvre} \begin{split} 
\br_1  &= \br_0 + \bnu_1 \times \bnu_0 , \\ \br_2 &= \br_0 - \bnu_2 \times \bnu_0 \end{split}
\end{equation} on every quad, where $\bnu_\ell = (-K(\br_\ell))^{-1/4}\bN_\ell$ for $\ell = 0,1,2,12$. These are the discrete Lelieuvre formula with non-constant curvature. In the case of constant negative curvature, these reduce to the classical discrete Lelieuvre equations first studied by Sauer \cite{sauer} and Wunderlich \cite{wunderlich}. In any case, \eqref{eq:discreteLelieuvre} guarantees that for all $(i,j)$, $\br_{i\pm1,j} - \br_{ij}$ and $\br_{i,j\pm1}-\br_{ij}$ are normal to $\bN_{ij}$, which is a discrete version of the geometry described by \eqref{eq:lelieuvre}.  


These equations can be enforced on any asymptotic complex by zooming into an elemetary quadrilateral and labeling one of the vertices as $\br_0$. This vertex has two neighbors: one along a $u$-edge which is labeled $\br_1$, and one along a $v$-edge which is labeled $\br_2$. The vertex diagonally opposite from $\br_0$ is labeled $\br_{12}$. 
Using \eqref{eq:discreteLelieuvre} we could arrive at two different formulas for $\br_{12}$ by either integrating from $0\to 1 \to 12$ or $0 \to 2 \to 12$. The discrete compatibility condition requires that these two formulas are equal. Following these calculations through, we see \begin{equation} \label{eq:r121}
\br_{12} =  \br_1 - \bnu_{12}\times \bnu_1 = \br_0 + \bnu_1 \times \bnu_0 - \bnu_{12} \times \bnu_1 = \br_0 -(\bnu_{12} + \bnu_0) \times \bnu_1.
\end{equation} and \begin{equation} \label{eq:r122}
\br_{12} =  \br_2 + \bnu_{12}\times \bnu_2 = \br_0 - \bnu_2 \times \bnu_0 + \bnu_{12} \times \bnu_2 = \br_0 + (\bnu_{12} + \bnu_0) \times \bnu_2.
\end{equation} Equating \eqref{eq:r121} and \eqref{eq:r122}, we arrive at the discrete compatibility condition: \begin{equation} \label{eq:discreteCompatibility}
(\bnu_{12} + \bnu_{0}) \times (\bnu_{1} + \bnu_{2}) = 0.
\end{equation} 

Interestingly, we can arrive at these exact same discrete equations and compatibility condition by faithfully discretizing \eqref{eq:lelieuvre} and \eqref{eq:compatibility}, which demonstrates a nice feature of discrete differential geometry: the discrete equations are not only analogs of the continuous equations, but actually establish an independent theory, which works in parallel with the continuous theory \cite{BobenkoSuris}. 

In the constant curvature case, \eqref{eq:discreteLelieuvre} and \eqref{eq:discreteCompatibility} reduce to \begin{equation} \label{eq:ddgConstCurv}\begin{split}
\br_1 = \br_0 &+ \bN_1 \times \bN_0, \,\,\,\,\,\,\,\,\,  \br_2 = \br_0 - \bN_2 \times \bN_0, \\
&(\bN_{12} + \bN_0) \times (\bN_1 + \bN_2) = 0.
\end{split}
\end{equation}
If we suppose that $\br_\ell$ and $\bN_\ell$ are known for $\ell = 0,1,2$,  then $\br_{12}$ and $\bN_{12}$ are uniquely and explicitly determined by \eqref{eq:ddgConstCurv}. Thus, upon appropriately assigning boundary data, one can use these equations to generate a discrete surface (in the form of an asymptotic complex) with constant negative curvature, as demonstrated by \cite{BobenkoPinkall,Toby,Yamamoto}. 

There is a wrinkle that must be addressed in the case of non-constant negative curvature: equations \eqref{eq:discreteLelieuvre} and \eqref{eq:discreteCompatibility} no longer explicitly determine $\br_{12}$ and $\bN_{12}$ in terms of known data. Indeed, suppose that for some elementary quadrilateral, $\br_\ell$ and $\bN_\ell$ (and thus $\bnu_\ell$) are known for $\ell = 0,1,2$. Our goal is to resolve $\br_{12}$ and $\bN_{12}$ using these values. To simplify notation, we introduce a radius of curvature $\rho(p) = (-K(p))^{-1/2}$, so that $\bnu_\ell = \rho_\ell^{1/2} \bN_\ell.$ From \eqref{eq:discreteCompatibility}, we see that \begin{equation} \label{eq:nuPreUpdate}
\bnu_{12} = \beta(\bnu_1 + \bnu_2) - \bnu_0. 
\end{equation} \\ To solve for $\beta$, we use $\|\bnu_\ell\|^2 =  \rho_\ell$. Thus \begin{equation} \rho_{12} = \beta^2 \|\bnu_1 + \bnu_2\|^2 - 2\beta \innerprod{\bnu_1 + \bnu_2}{\bnu_0} + \rho_0, \end{equation} so that \begin{equation} \label{eq:betaFormula} \begin{split} \beta &= \Big(1 + \sqrt{1+\alpha}\Big) \frac{\innerprod{\bnu_1 + \bnu_2}{\bnu_0}}{\|\bnu_1 + \bnu_2\|^2}, \\ \alpha &= \frac{\|\bnu_1 + \bnu_2\|^2 (\rho_{12} - \rho_0)}{\innerprod{\bnu_1 + \bnu_2}{\bnu_0}^2},\end{split} \end{equation} and finally \begin{equation} \label{eq:normalUpdateEquation}
\bN_{12} = \sqrt{\frac{\rho_0}{\rho_{12}} }\bigg( \Big(1+\sqrt{1+\alpha}\Big) P - I \bigg)\bN_0
\end{equation} where \begin{equation}\label{eq:P}
P = \frac{(\bnu_1 + \bnu_2)(\bnu_1 + \bnu_2)^t}{\|\bnu_1 + \bnu_2\|^2}.
\end{equation} Note that $P$ is rank 1, and so \eqref{eq:normalUpdateEquation} expresses $\bN_{12}$ as perturbation of a Householder reflection of $\bN_0$ about the plane with normal $\bnu_1 + \bnu_2$. The perturbations are expressed by the quotient $(\rho_{0}/\rho_{12})^{1/2}$ and the parameter $\alpha \propto (\rho_{12} - \rho_0)$. In particular, if the curvature is constant, both of these perturbations disappear, and this is a standard Householder reflection.

To summarize, if $\br_{\ell}$ and $\bN_\ell$ are known for $\ell = 0,1,2$, we can resolve $\bN_{12}$ and $\br_{12}$ using \begin{equation} \label{eq:fullUpdate}
\begin{split}
\bN_{12} &= \sqrt{\frac{\rho_0}{\rho_{12}} }\bigg( \Big(1+\sqrt{1+\alpha}\Big) P - I \bigg)\bN_0, \\
\br_{12} &= \br_{2} + \bnu_{12} \times \bnu_{2} \,\,\,\,\,\,\,\,\,\,\,\,\, (\text{or } \,\,\,\,\,\,\br_{12} = \br_1 - \bnu_{12} \times \bnu_{1}).
\end{split}
\end{equation} Once boundary data is assigned, one can use this system to ``step inward'' from the boundary and generate an asymptotic complex with prescribed Gaussian curvature. Note that this is an implicit system of equations: $\br_{12}$ depends on $\bN_{12}$ through $\bnu_{12}$, while $\bN_{12}$ depends on $\br_{12}$ through $\alpha$ and $\rho_{12}$. One way around this would be to let the curvature $K$ depend directly on the asymptotic coordinates: $K(u,v)$ rather than $K(p)$. In this case, the equations become explicit. However, this may be unsatisfying since ultimately the curvature should be intrinsic, while the asymptotic coordinates depend on the particular immersion into $\R^3$. 

Alternatively, we propose solving \eqref{eq:fullUpdate} by iteration. This could be done in a few different ways. For example, one could solve \eqref{eq:fullUpdate} by fixed point iteration at each quadrilateral as one generates the surface. By constrast, one could generate an entire discrete surface and normal field $(\br_{ij}, \bN_{ij})$ using a fixed (and known) curvature function, reset the Gaussian curvature $K_{ij} = K(\br_{ij})$ and use this curvature to regenerate the surface.  We opt for the latter method. In this manner, we are only ever solving explicit equations since $\rho_{12}$ and $\alpha$ will be treated as known whenever they are encountered. Thus, rather than approximately solving \eqref{eq:fullUpdate} on each quadrilateral, we exactly solve \eqref{eq:fullUpdate} but with approximations of $\alpha$ and $\rho_{12}$.

We are interested in a situation where curvature increases as the surface expands outward from some central point, which we think of as the ``origin'' of the surface.  In this case, curvature is not an explicit function of the coordinate $\br$ but rather of the geodesic distance from $\br$ to the origin. This  is inspired by many physical examples of negatively curved surfaces such as leafy lettuce, marine flatworms, and flowers like the red bird of paradise pictured in figure \ref{fig:rbop}. In these cases and many others, the surface is nearly flat toward the middle, but becomes curved and wrinkled near the edges. 

We generate a surface whose curvature increases with geodesic distance from the center by viewing the surface as a perturbation of a pseudospherical surface. Specifically, supposing we already have a discrete surface with Gaussian curvature $K(p)$, that $D(p)$ is the geodesic distance from the center of the surface to the point $p$ and $h:[0,\infty) \to [0,\infty)$ is a nondecreasing perturbation function, Algorithm \ref{alg:IterationAlg} gives a manner of generating a discrete surface with Gaussian curvature \begin{equation}\label{eq:curvatureEps} K_{\eps}(p) = K(p) - \eps h(D(p)).\end{equation}  

\begin{algorithm}[t!]
\caption{Generating a discrete surface with prescribed negative curvature $K_\eps(p)$ given by \eqref{eq:curvatureEps}, assuming a discrete surface with Gaussian curvature $K(p)$ is already resolved.}

Specify $\eps > 0$ and $h:[0,\infty)\to[0,\infty)$ Lipschitz continuous to define $K_\eps(p)$ as in \eqref{eq:curvatureEps}. Set $\text{CHANGE} = \infty$, specify a convergence tolerance $\text{TOL} > 0$, and set $n = 0$.\\

Also input $\br^{(0)}_{ij}$ and $\bN^{(0)}_{ij}$: a discrete immersion and normal field with Gaussian curvature $K(p)$, and let $K_{ij}$ represent the curvature at $\br^{(0)}_{ij}$. \\

\begin{algorithmic}
\While {$\text{CHANGE} > \text{TOL}$}\\
    \State Compute the geodesic distance function $D^{(n)}_{ij}$ which gives geodesic distance from $\br^{(n)}_{ij}$ \\ \hspace*{\algorithmicindent}to the origin.\\

    \State Assign $K^{(n)}_{ij} \longleftarrow K_{ij}-\eps h(D^{(n)}_{ij})$\\

    \State Generate a new surface $\br^{(n+1)}_{ij}$ and normal field $\bN^{(n+1)}_{ij}$ using \eqref{eq:fullUpdate} with curvature $K^{(n)}_{ij}$\\ \hspace*{\algorithmicindent}and the same boundary data as used for initialization.\\

    \State Assign $\text{CHANGE} \longleftarrow \| \br^{(n)}_{ij} - \br^{(n+1)}_{ij}\|$\\

    \State Assign $n \longleftarrow n+1$\\
\EndWhile

\State{\bf return} $\br_{ij} = \br^{(n+1)}_{ij}, \bN_{ij} = \bN^{(n+1)}_{ij}$. These are the coordinates and normal field for a discrete surface with approximate curvature given by $K_{\eps}(p)$.
\end{algorithmic}
\label{alg:IterationAlg}
\end{algorithm}

The convergence criterion used in the algorithm is somewhat arbitrary. For our purposes, we halt the iteration when the maximal difference between successive surfaces $\max_{ij} \|r^{(n+1)}_{ij} - r^{(n)}_{ij} \|$ falls beneath some tolerance. We will prove that, under some additional conditions specified below, the iteration in Algorithm \ref{alg:IterationAlg} will converge for the Amsler-type surfaces of interest in Section \ref{sec:amsler}, given that the perturbation size $\eps$ is small enough. In practice, we will always begin with a pseudospherical surface which can be generated using \eqref{eq:ddgConstCurv}, once boundary data has been initialized as specified below. There is no reason one need only consider small perturbations from pseudospherical surfaces. To generate a surface with arbitrary (Lipschitz continuous) curvature function $K\le -1$, we first generate a pseudospherical surface, and then write $$K = -1 - (-1-K) \backdefeq -1 - h = -1 -\sum_{n=1}^N \frac h N.$$ Algorithm \ref{alg:IterationAlg} then allows to successively generate surfaces with discrete curvature given by $$-1, \,\, -1-\frac{h}N, \,\, -1-\frac{2h}N,\,\, -1-\frac{3h}N,\,\, \ldots$$ by successively applying Algorithm \ref{alg:IterationAlg} with $K_{ij} = -1-nh_{ij}/N$ and $\eps = 1/N$.  Because the proof of convergence for Algorithm \ref{alg:IterationAlg} is somewhat long and technical, we relegate it to Appendix \ref{sec:convProof}. 

\section{Geodesic Distance on Triangulated Manifolds}\label{sec:geodesic}

Lastly, to run the above algorithm and generate surfaces of non-constant negative curvature, we need to compute the geodesic distance function on the surface.  Accordingly, we present a method for finding the geodesic distance on a triangulated 2-manifold $M$ embedded in $\mathbb R^3$. Given a source point $o \in M$ (which we think of as the ``origin'' of the manifold), the geodesic distance function on the manifold solves the Eikonal equation $$\abs{\nabla D(x)} = 1, \,\,\,\,\ x \in M \setminus \{o\},$$  with the condition that $D(o) = 0$. A classic approach due to Kimmel and Sethian \cite{KimmelSethian} is to directly approximate this equation on a grid and sweep through the nodes in the fast marching manner \cite{Sethian1,SethVlad1,Tsitsiklis}. We modify their algorithm, maintaining the fast marching methodology, but rather than approximating the Eikonal equation, we imagine ``unfolding'' the manifold and explicitly calculating the geodesic distance to nodes using the known values of the geodesic distance at nearby nodes. 

The idea behind the fast marching method is to track information as a wavefront propagating outward from boundary data and update nodes as the wavefront reaches them. Explicitly, suppose we have a discrete triangulation $\{ r_i\}_{i \in I}$ of our manifold embedded in $\R^3$, and suppose our triangulation contains the origin $o$. Let $D_i \approx D(r_i)$ denote our approximation of the distance function at the node $r_i$. The basic fast marching algorithm is contained in algorithm \ref{alg:fastMarch}.

\begin{algorithm}[t!]
\caption{The basic fast marching methodology for geodesic distance.}

Begin with a discrete triangulation $\{r_i\}_{i\in I}$ of a manifold embedded in $\R^3$, and specify the node $r_o$ which will serve as the ``origin."  Set $D_o = 0$ and mark $r_o$ as \textit{Accepted}. Mark all other nodes $r_i$ as \textit{Far Away} and set $D_i = +\infty$. \\

\begin{algorithmic}
\While {there are still nodes which are not \textit{Accepted}} \\

\State Mark any node $r_i$ which has an \textit{Accepted} neighbor as \textit{Considered}. \\

\State At each \textit{Considered} node, update the approximation $D_i$. (We describe our update \\ \hspace*{\algorithmicindent}procedure below).\\

\State Find the \textit{Considered} node $r_i$ with the smallest approximate distance $D_i$. Remove this \\\hspace*{\algorithmicindent}node from the \textit{Considered} list, and mark it as \textit{Accepted}.\\

\EndWhile

\State{\bf return} the approximate discrete geodesic distance function $D_i$. 
\end{algorithmic}
\label{alg:fastMarch}
\end{algorithm}

If there are $N$ nodes, the algorithm has complexity $O(N\log N)$, where the $\log N$ represents the effort to maintain a list of \emph{Considered} nodes which is sorted by approximate value $D_i$ \cite{SethVlad1}.

To use the algorithm, one needs to specify an update rule for obtaining $D_i$ from the values of neighboring nodes. As stated above, this can be done by discretizing the Eikonal equation as in \cite{KimmelSethian, Sethian1}. Alternatively, one can exploit the geometry and use the following procedure. 

\begin{figure}[b!]
\centering
\includegraphics[width = 0.6\textwidth,trim= 80 100 60 50,clip]{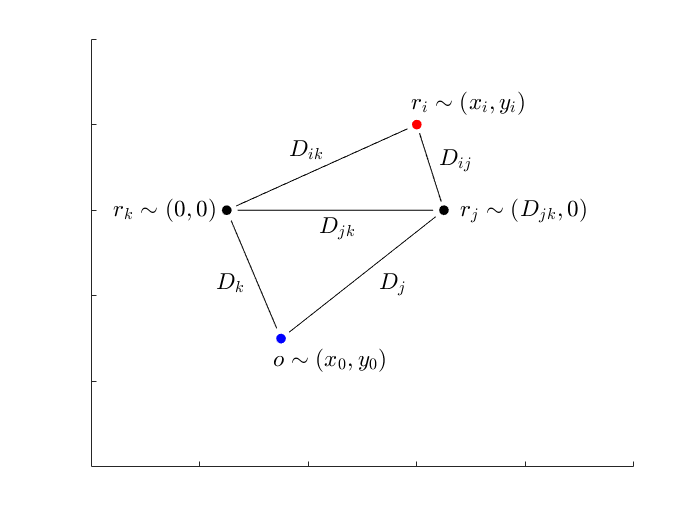}
\caption{``Unfolding'' the manifold and locating different nodes in the $xy$-plane.}
\label{fig:plane}
\end{figure}

Suppose the node $r_i$ has neighbors $r_j, r_k$ whose geodesic distances $D_j, D_k$ are known. For any indices $i,j$ let $D_{ij} = \|r_i - r_j\|$. Consider an ``unfolding'' and reorienting of the manifold which places it over the typical $xy$-Cartesian plane with $r_k$ at $(0,0)$ and $r_j$ on the $x$-axis at $(D_{jk},0)$. The goal is to explicitly locate $o$ and $r_i$ in this plane, using the known distances $D_j, D_k, D_{ij}, D_{ik}, D_{jk}$. This is pictured figure \ref{fig:plane}.

Under this unfolding, the coordinates for the origin satisfy the equations 
\begin{equation}\label{eq:originSystem} 
\begin{split}
 x_o^2 + y_o^2 &= D_k^2,\\ 
 (x_o - D_{jk})^2 + y_o^2 &= D_j^2.
 \end{split}
 \end{equation}
  Likewise, the coordinates for the node $r_i$ satisfy 
  \begin{equation}\label{eq:riSystem} 
\begin{split}
 x_i^2 + y_i^2 &= D_{ik}^2,\\ 
 (x_i - D_{jk})^2 + y_i^2 &= D_{ij}^2.
 \end{split}
 \end{equation} 
 Solving these quadratic systems yields 
 \begin{equation}\label{eq:originCoords} 
\begin{split}
 x_o &= \frac{D_k^2 - D_j^2 + D_{jk}^2}{2D_{jk}},\\ 
y_o^{\pm}&= \pm \frac{1}{2D_{jk}} \sqrt{2D_{jk}^2 D_k^2 - D_{jk}^4 + 2D_{jk}^2D_j^2 - (D_j^2 - D_k^2)^2}.
 \end{split}
 \end{equation} and
  \begin{equation}\label{eq:riCoords} 
\begin{split}
 x_i &= \frac{D_{ik}^2 - D_{ij}^2 + D_{jk}^2}{2D_{jk}},\\ 
y_i^{\pm} &= \pm \frac{1}{2D_{jk}} \sqrt{2D_{jk}^2 D_{ik}^2 - D_{jk}^4 + 2D_{jk}^2D_{ij}^2 - (D_{ij}^2 - D_{ik}^2)^2}.
 \end{split}
 \end{equation} 
This gives different possibilities for where the $o$ and $r_i$ lie with respect to $r_j$ and $r_k$, and thus different approximations $$D_i = \sqrt{(x_i - x_o)^2 + (y_i^{\pm} - y_o^{\pm})^2}.$$ However, we only need to use the configuration that suggests the largest distance from $o$ to $r_i$. The reasoning for this is that characteristics of the Eikonal equation flow outward from the origin. If $r_i, r_j, r_k$ are as pictured in figure \ref{fig:plane}, then $D_i \ge \max\{D_j, D_k\}$, and the maximum of the suggested distances $$D^{(jk)}_i = \sqrt{(x_i - x_o)^2 + (y_i^+ - y_o^-)^2}$$ is the correct approximation to $D(r_i)$ from the triangle $(r_i,r_j,r_k)$. If instead $r_i$ was closer to the origin than one of $r_j, r_k$, then this formula provides an incorrect approximation---one which is too large. But in that case, information from the farther node should not influence our approximation of $D_i$. Accordingly, we check all triangles $(r_i,r_j,r_k)$ in our mesh containing $r_i$, find the suggested approximation $D^{(jk)}_{i}$ from this triangle, and take the minimum of all suggested approximations. 

There is one final way in which the picture in figure \ref{fig:plane} may be faulty. We could also have a situation where both $r_j, r_k$ are closer to the origin than $r_i$, but the line from $r_i$ to the origin (and thus the characteristic that determines the value of $D_i$) does not pass through the triangle $(r_i,r_j,r_k)$. In this case, the appropriate approximation of $D_i$ from the triangle $(r_i,r_j,r_k)$ should not be given by the formula above, but rather by traveling directly from $r_j$ or $r_k$ to $r_i$; that is, $D^{(jk)}_i = \min\{D_j + D_{ij}, D_k + D_{ik}\}.$ Practically, this update will be used at the beginning, when the only known information is $D_o = 0$. It will also be used whenever $r_j$ is the closest node to $r_i$ and $o,r_i,r_j$ are colinear. 

Putting this all together, our approximation to $D_i$ from within the triangle $(r_i,r_j,r_k)$ is \begin{equation}\label{eq:distanceUpdate} D^{jk}_i = \min\bigg\{ \sqrt{(x_i - x_o)^2 + (y_i^+ - y_o^-)^2}, D_j + D_{ij}, D_{k} + D_{ik}\bigg \}\end{equation} where $(x_o,y_o)$ and $(x_i,y_i)$ are given in \eqref{eq:originCoords} and \eqref{eq:riCoords} respectively. The final approximation of $D_i$ is then $$D_i = \min_{j,k} D^{(jk)}_i,$$ where the minimum is taken with respect to all pairs of indices $(j,k)$ such that $(r_i,r_j,r_k)$ is a triangle in our mesh and the distances $D_j,D_k$ have already been \emph{Accepted}.

One final note is that there is no reason the ``origin" needs to be a single point. It could rather be a set of points where the distance function is known, and the algorithm works identically except that during the initialization, all ``origin" points should be immediately marked as \textit{Accepted} and their distances should be fixed.

\begin{figure}[t!]
\centering
\includegraphics[width = 0.4\textwidth]{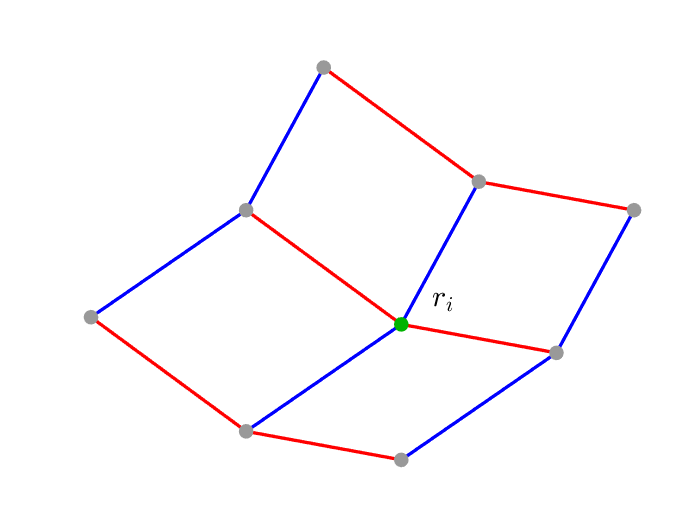}\,\, \includegraphics[width = 0.4\textwidth]{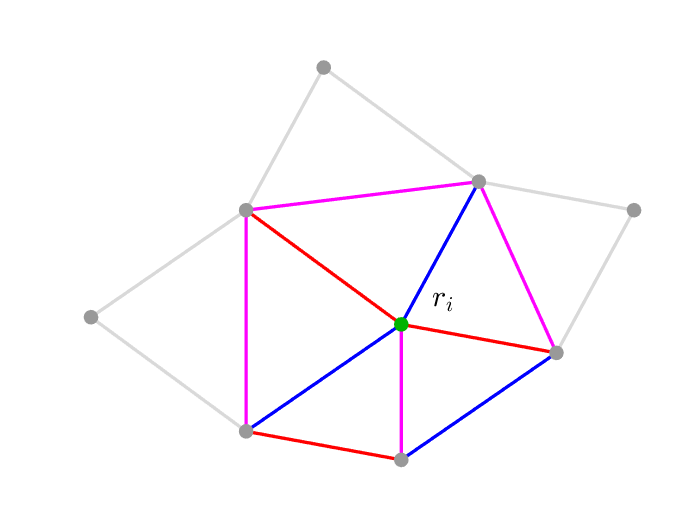} 
\caption{Left: asymptotic complex near a node $r_i$. Right: the local triangulation resulting from bisecting obtuse angles. Blue and red lines are $u$- and $v$-asymptotic lines respectively. Magenta lines are artificial lines used for triangulation.  The gray lines (right) belong to quads which abut $r_i$ in the asymptotic complex, but do not belong to triangles which abut $r_i$ in the triangulation.}
\label{fig:triangulation}
\end{figure}

To use this method for finding the distance function, we need to triangulate our asymptotic complexes. This can be done by drawing artificial lines along the $(u+v)$-direction or $(u-v)$-direction, where $u$ and $v$ are the asymptotic coordinates. Convergence for fast marching schemes is only guaranteed when the triangulation of the manifold is acute. Accordingly, when we triangulate our manifold, we check which angles need to be split so that the triangulation is indeed acute. This is shown in figure \ref{fig:triangulation}. As pictured, the blue and red lines are $u$- and $v$-asymptotic lines respectively, and the magenta lines are artificial lines used to triangulate the mesh. Note that obtuse angles are bisected so that the resulting triangulation is acute. The gray edges in the right of figure \ref{fig:triangulation} belong to quads for which $r_i$ is a vertex in the asymptotic complex, but do not belong to triangles for which $r_i$ is a vertex in the triangulation.   \\

\section{DDG for Amsler-Type Surfaces} \label{sec:amsler}

We demonstrate our methods using Amsler-type surfaces \cite{Amsler}. Such surfaces are treated in detail in Chap.~5 of Ref.~\cite{BobenkoUlrich}. Amsler surfaces are comprised of sectors, each of which can be seen as an immersion $\br: [0,\infty) \times [0,\infty)\to \R^3$ which obeys the Lelieuvre formulas, as well as the condition that the bounding asymptotic curves for the sector, which are given by $\br(0,\cdot)$ and $\br(\cdot,0)$, are geodesics {\em in the ambient space} $\mathbb{R}^3$. Consequently, these curves are necessarily also geodesics on the surface with respect to its intrinsic geometry. As discussed by \cite{Toby}, one can generate such surfaces sector-by-sector and, if the angle between the asymptotic lines is chosen appropriately, paste the sectors together such that the normal $\bN$ remains continuous across boundary lines.

We form discrete Amsler sectors by instead considering immersions $\br :\mathbb N^2 \to \R^3$. In this case, asymptotic curves are the images of $\{(i,j_0)\}_{i\in\mathbb N}$ and $\{(i_0,j)\}_{j\in\mathbb N}$ under the immersion $\br$. One views a single sector as the image of the rectangular grid $\{(i,j) \, : \, 0 \le i \le I, 0 \le j \le J\} \subset \mathbb Z^2$. The two necessary pieces of boundary data are (1) the angle $\phi_1$ between the asymptotic lines at the origin and (2) the units of distance in the $u$- and $v$-directions. We describe the procedure here.

\begin{figure}[t!]
\centering
\begin{subfigure}{0.4\textwidth}
\centering
\includegraphics[width=\textwidth]{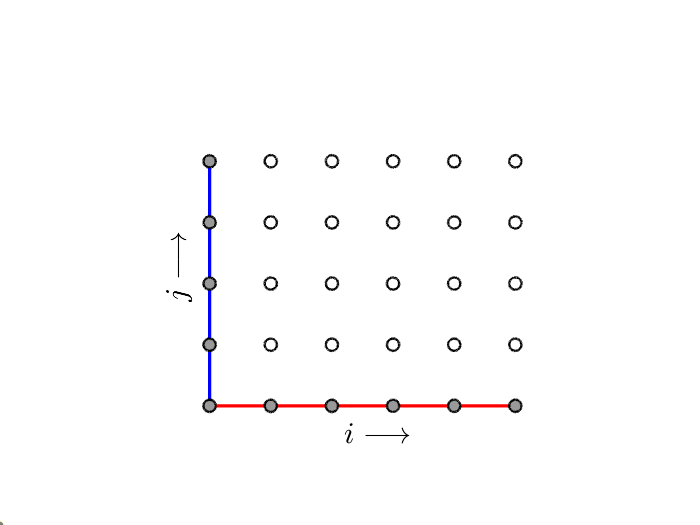}
\caption{Domain for the immersion $\br$. Filled circles indicate positions where data is provided.}
\label{fig:ijgrid}
\end{subfigure}\hfill
\begin{subfigure}{0.4\textwidth}
\centering
\includegraphics[width=\textwidth]{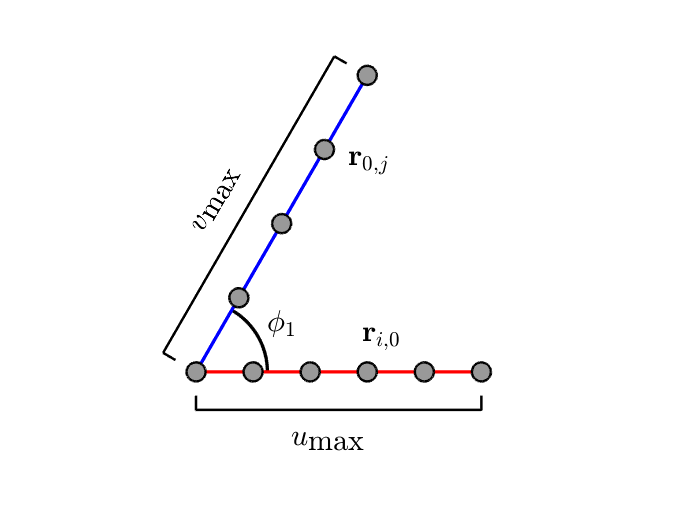}
\caption{Boundary of the Amsler-type sector. That is, image of the immersion $\br$, plotted in the $xy$-plane.}
\label{fig:uvgrid}
\end{subfigure}
\caption{Our immersion $\br$ will map a rectangle in $\Z^2$ into a sector of angle $\phi_1$ in the $xy$-plane. The boundary coordinates $\br_{i,0}$ and $\br_{0,j}$ lie in the $xy$-plane. In general, $\br_{ij}$ will have a nonzero $z$-coordinate.}
\label{fig:amslerSurfaceBCs}
\end{figure}

Given the angle $\phi_1$ between the asymptotic lines and the maximal distances $u_{\text{max}}$ and $v_{\text{max}}$ in the $u$- and $v$-directions, the boundary vertices $\{\br_{i,0}\}_{0\le i \le I}$ and $\{\br_{0,j}\}_{0\le j \le J}$ of the ``principal'' sector will be comprised (respectively) of evenly spaced points along the $x$-axis from the origin to distance $u_{\text{max}}$ and along the line $y/\sin \phi_1 = x/\cos \phi_1$ from the origin to distance $v_{\text{max}}$. That is, the images of the boundary lines remain in the $xy$-plane. We initialize by specifying lengths rather than angles so that the geodesic distance function along these boundary lines is explicitly known. This is pictured in figure \ref{fig:amslerSurfaceBCs}. In the constant curvature case, with boundary data for the normal field appropriately chosen, the quadrilaterals comprising the surface are parallelograms and the normal vectors at the vertices of the quadrilaterals form spherical parallelograms \cite{Toby}. Thus there are some angles $\delta_u, \delta_v$ such that \begin{equation} \label{eq:BCsAngle} \|\br_{i\pm1,0} - \br_{i,0}\| = \|\bN_{i\pm 1,0}\times \bN_{i,0}\| = \sin \delta_u \,\,\,\,\, \text{ and } \,\,\,\,\, \| \br_{0,j\pm1}-\br_{0,j}\| = \|\bN_{0,j \pm 1}\times \bN_{0,j}\| = \sin \delta_v, \end{equation} and it follows that $\innerprod{\bN_{i\pm1,0}}{ \bN_{i,0}} = \cos \delta_u, \,\, \innerprod{\bN_{0,j\pm1}}{\bN_{0,j}} = \cos \delta_v$. We choose the boundary data for our normal field to mimic this behavior: that is, we fix the angles between the successive normal vectors on along the boundary. Specifically, we set $\bN_{0,0} = (0, 0, 1)$ so that near the origin the surface looks approximately like the $xy$-plane. We then progressively rotate $\bN_{i,0}$ about the $x$-axis, and progressively rotate $\bN_{0,j}$ about the line $y/\sin \phi_1 = x/\cos \phi_1$, maintaining \begin{equation} \label{eq:normalBCs1}
\innerprod{\bN_{i\pm1,j}}{ \bN_{ij}} = \cos \delta_u \,\,\,\,\,  \text{ and } \,\,\,\,\, \innerprod{\bN_{i,j\pm1}}{\bN_{ij}} = \cos \delta_v,
\end{equation} where $\delta_u,\delta_v$ are as in \eqref{eq:BCsAngle}. However, in our case, we have \begin{equation} \|\br_{i\pm 1,j} - \br_{i,j}\| = (\rho_{i\pm 1,j}\rho_{i,j})^{1/2} \|\bN_{i\pm 1,j} \times \bN_{i,j}\|\end{equation}  and similarly for $\|\br_{i, j\pm 1} - \br_{i,j}\|$.  In this case, the angles $\delta_u$ and $\delta_v$ are not constant, but depend on the prescribed curvature (which depends on geodesic distance). Sorting out the algebra, we initialize the boundaries of normal field using \begin{equation}\label{eq:normalBCs}
\begin{split}
\bN_{i+1,0} &= \text{Rot}(\delta_{i} \,;\, (1,0,0))\bN_{i,0}, \\ 
\bN_{0,j+1} &= \text{Rot}(\delta_j\,;\, -(\cos(\phi_1), \sin(\phi_1),0) )\bN_{0,j},
\end{split}
\end{equation} where $\text{Rot}(\delta,\mathbf{a})$ denotes the rotation by angle $\delta$ about the vector $\mathbf{a}$, and \begin{equation}\label{eq:rotationAngles}
\begin{split}
\delta_i &= \arcsin\left( \frac{\|\br_{i+1,0} - \br_{i,0}\|}{(\rho_{i+1,0}\rho_{i,0})^{1/2} }\right), \\ 
\delta_j &= \arcsin\left( \frac{\|\br_{0,j+1} - \br_{0,j}\|}{(\rho_{0,j+1}\rho_{0,j})^{1/2} }\right)
\end{split}
\end{equation} Because we are specifying the distances $\|\br_{i+1,0} - \br_{i,0}\|$ and $\|\br_{0,j+1} - \br_{0,j+1}\|$, $\rho_{i,0}$ and $\rho_{0,j}$ are known, so these formulas are entirely explicit, and it is straightforward to verify that these definitions do indeed satisfy the discrete Lelieuvre equations \eqref{eq:discreteLelieuvre} along the boundary lines.

Using this initialization, we generate and Amsler sector by iterating \eqref{eq:fullUpdate} with prescribed curvature \eqref{eq:curvatureEps} as detailed in algorithm \ref{alg:IterationAlg}. As described, this generates an Amsler sector whose boundary lines extend in the directions of $\B s_1 = (1,0,0)$ and $\B s_2 = (\cos(\phi_1),\sin(\phi_1),0)$. In an analogous manner, one can generate another Amsler surface whose boundary lines extend in the directions $\B s_2 = (\cos(\phi_1),\sin(\phi_1),0)$ and $\B s_3 = (\cos(\phi_2),\sin(\phi_2),0)$ for some $\phi_2 > \phi_1$, and if the boundary data is consistent along the mutual boundary line $\B s_2$, these sectors can be patched together so that the normal maps remain continuous across the boundary. 

\begin{figure}[b!]
\centering
\includegraphics[width=0.5\textwidth]{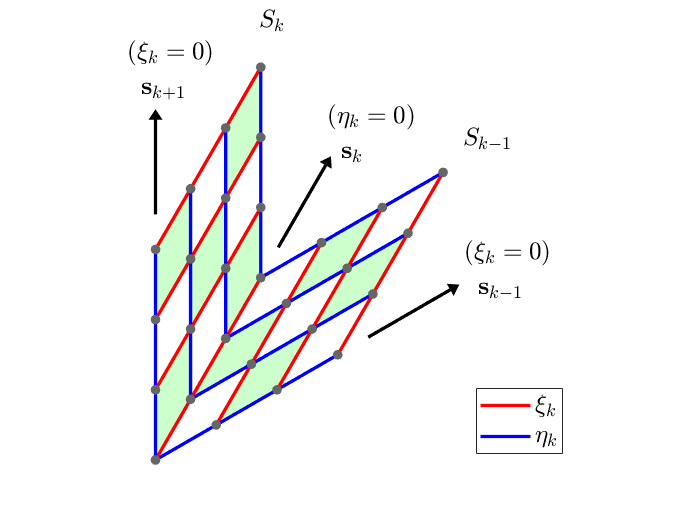}
\caption{Patching together of adjoining Amsler sectors (projected down to the plane). The coordinates $(\xi_k,\eta_k)$ describe the pair of adjoining sectors. The origin is at $(\xi_k,\eta_k) = (0,0)$, the line in the direction of $\B s_k$ is the line $\eta_k = 0$ and the V-shape formed by adjoining $\B s_{k+1}$ and $\B s_{k-1}$ is the line $\xi_k = 0$, with $\B s_{k+1}$ being captured by $\eta >0$ and $\B s_{k-1}$ by $\eta_k < 0$. }
\label{fig:patching}
\end{figure}

This sort of patching is described in detail in \cite[Sect. 3.2]{Toby}. Generally, for any angles $\theta_k \in (0,\pi)$ for $k = 1,\ldots,2n \ge 4$ such that $\sum_{k=1}^{2n} \theta_k = 2\pi$, one can generate a collections of Amsler sectors $\{S_k\}_{k=1}^{2n}$ whose boundary lines extend from the origin in the directions of $\B s_k = (\cos(\phi_k),\sin(\phi_k),0)$, where $\phi_k = \sum^k_{\ell=1} \theta_\ell.$ The line in the direction of $\B s_k$ is a border shared by sectors $S_{k-1}$ and $S_{k}$ (with the obvious adjustment when $k=1$). One can then define coordinates $(\xi_k,\eta_k)$ on the union of adjoining sectors $S_{k-1} \cup S_k$, where $\xi_k \ge 0$ denotes distance along $\B s_k$ and $\abs{\eta_k}$ denotes distance along $\B s_{k+1}$ if $\eta_k > 0$ or distance along $\B s_{k-1}$ if $\eta_k < 0.$ This is demonstrated in figure \ref{fig:patching}, where we plot two adjoining Amsler sectors (projected down to the plane). In figure \ref{fig:patching}, the coordinates $(\xi_k,\eta_k)$ describe the pair of adjoining sectors. Red lines are curves of constant $\eta_k$, and the blue V-shapes are curves of constant $\xi_k$. 

Having defined these coordinates, the asymptotic coordinates on sector $S_k$ are given by $(u,v) = (\xi_k,\eta_k)$ if $k$ is odd, and $(u,v) = (\eta_k,\xi_k)$ if $k$ is even. Using these asymptotic coordinates and the boundary values described above, the Lelieuvre equations define immersions and normal fields for each sector. It is demonstrated in \cite{Toby} that the normal fields for adjoining sectors do indeed remain continuous across the boundaries as the sectors are pasted together. 

\begin{figure}[t!]
\begin{subfigure}{0.31\textwidth}
\centering 
\includegraphics[width=\textwidth,trim = 80 80 80 80, clip]{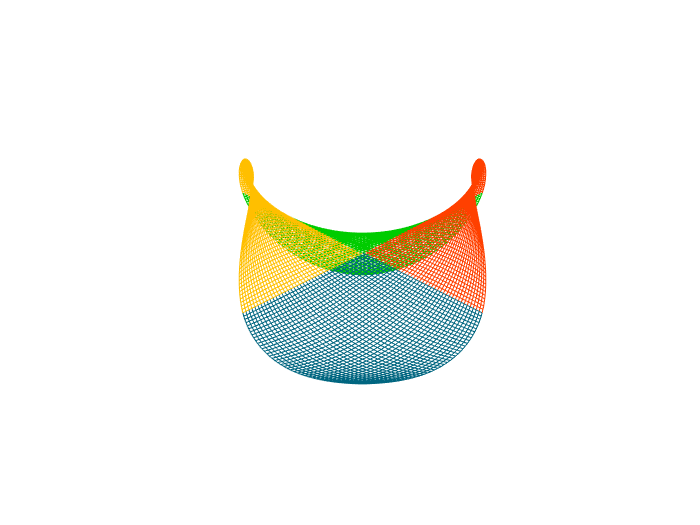}
\caption{$n=2, \eps = 1$}
\label{fig:n2eps1}
\end{subfigure}~
\begin{subfigure}{0.31\textwidth}
\centering 
\includegraphics[width=\textwidth,trim = 80 80 80 80, clip]{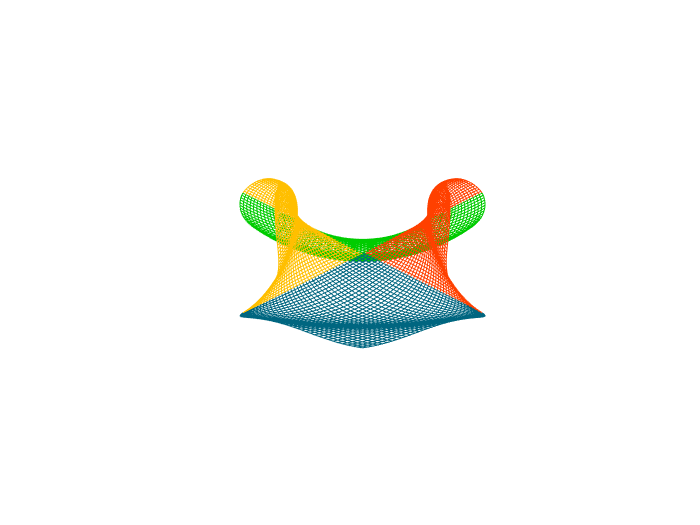}
\caption{$n=2, \eps = 10$}
\label{fig:n2eps10}
\end{subfigure}~
\begin{subfigure}{0.31\textwidth}
\centering 
\includegraphics[width=\textwidth,trim = 80 80 80 80, clip]{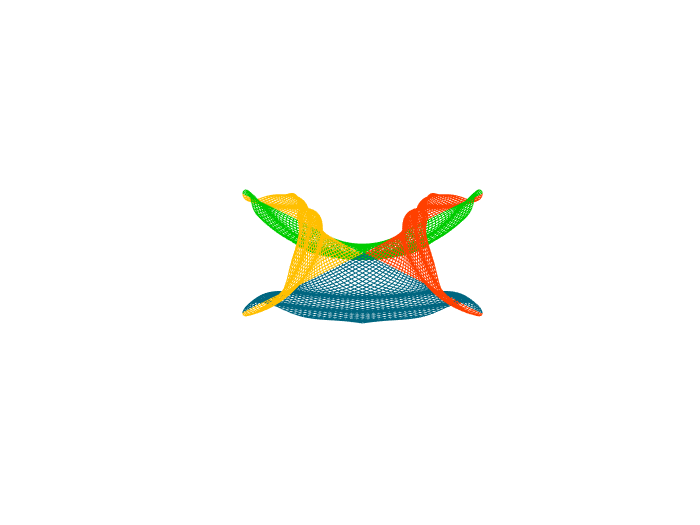}
\caption{$n=2, \eps = 50$}
\label{fig:n2eps50}
\end{subfigure}~ \\
\begin{subfigure}{0.31\textwidth}
\centering 
\includegraphics[width=\textwidth,trim = 80 40 80 40, clip]{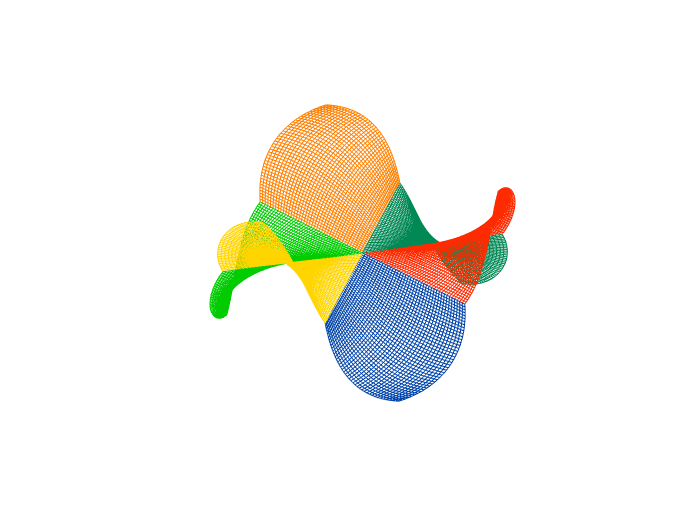}
\caption{$n=3, \eps = 1$}
\label{fig:n3eps1}
\end{subfigure}~
\begin{subfigure}{0.31\textwidth}
\centering 
\includegraphics[width=\textwidth,trim = 80 40 80 40, clip]{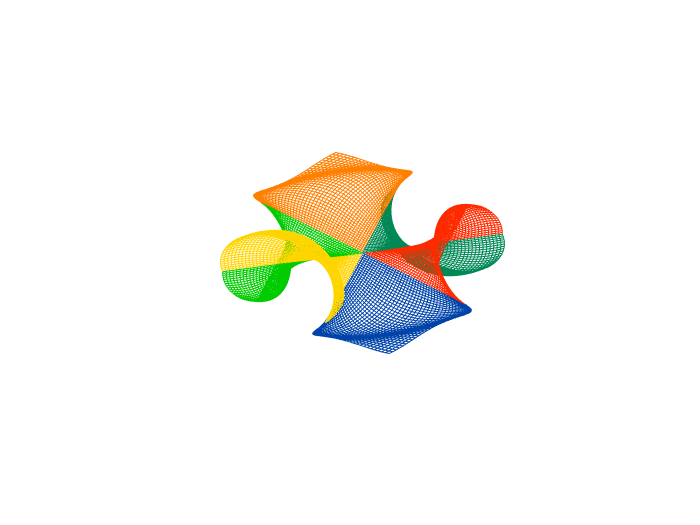}
\caption{$n=3, \eps = 10$}
\label{fig:n3eps10}
\end{subfigure}~
\begin{subfigure}{0.31\textwidth}
\centering 
\includegraphics[width=\textwidth,trim = 80 40 80 40, clip]{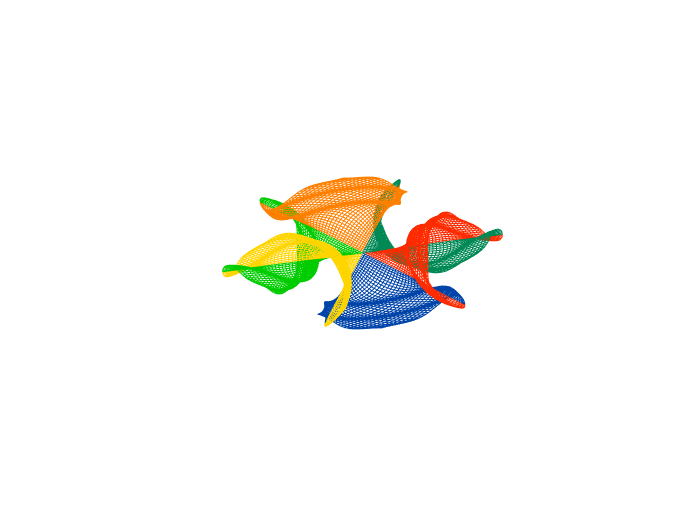}
\caption{$n=3, \eps = 50$}
\label{fig:n3eps50}
\end{subfigure}~ \\
\begin{subfigure}{0.31\textwidth}
\centering 
\includegraphics[width=\textwidth,trim = 80 40 80 40, clip]{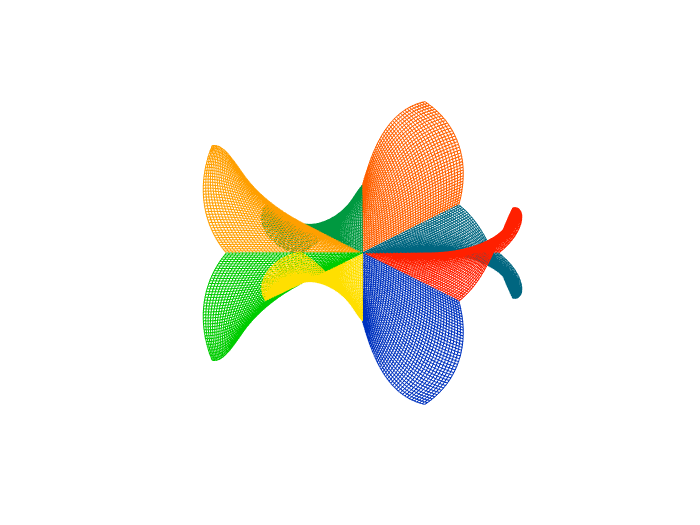}
\caption{$n=4, \eps = 1$}
\label{fig:n4eps1}
\end{subfigure}~
\begin{subfigure}{0.31\textwidth}
\centering 
\includegraphics[width=\textwidth,trim = 80 40 80 40, clip]{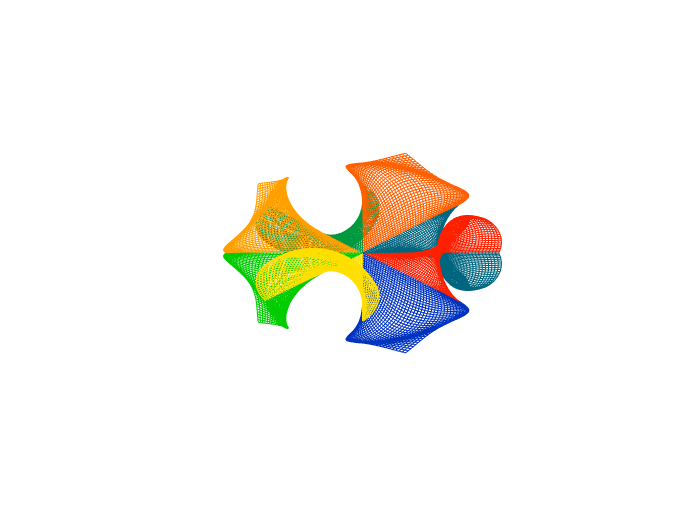}
\caption{$n=4, \eps = 10$}
\label{fig:n4eps10}
\end{subfigure}~
\begin{subfigure}{0.31\textwidth}
\centering 
\includegraphics[width=\textwidth,trim = 80 40 80 40, clip]{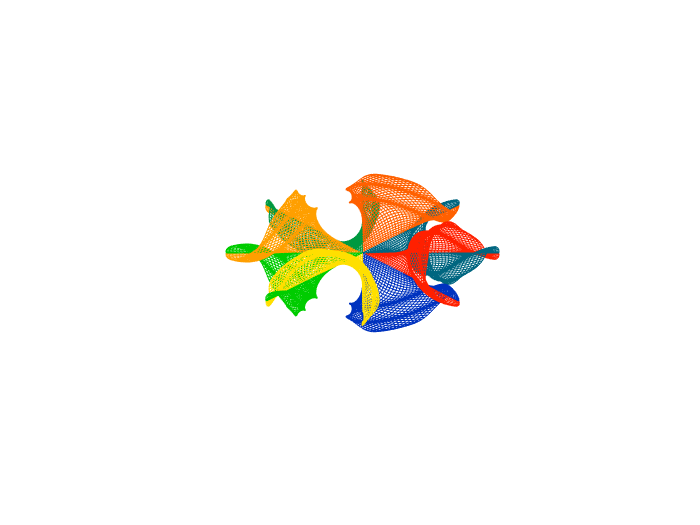}
\caption{$n=4, \eps = 50$}
\label{fig:n4eps50}
\end{subfigure}~ \\
\caption{Discrete Amsler surfaces with $2n$ sectors and prescribed curvature $K(p) = -(1+\eps D(p))$.}
\label{fig:amslerSurfaces1}
\end{figure}

Our figures demonstrate the maximally symmetric case that $\theta_k = \pi/n$ for all $k = 1,\ldots, 2n$. Examples of surfaces are included in figure \ref{fig:amslerSurfaces1}. Here we use $n=2,3$ and $4$ resulting in $4,6$ and $8$ sectors respectively. We prescribe the curvature to increase linearly with geodesic distance: $K_\eps(p) = -(1+\eps D(p))$ where $D(p)$ is geodesic distance to the origin. As we increase $\eps$, the resulting surface becomes increasingly curved, especially toward the boundary. Note, we generated these surfaces using the iterative method and stepping in $\eps$ as described above. 

\begin{figure}[t!]
\begin{subfigure}{0.31\textwidth}
\centering
\includegraphics[width=\textwidth,trim = 120 120 120 120, clip]{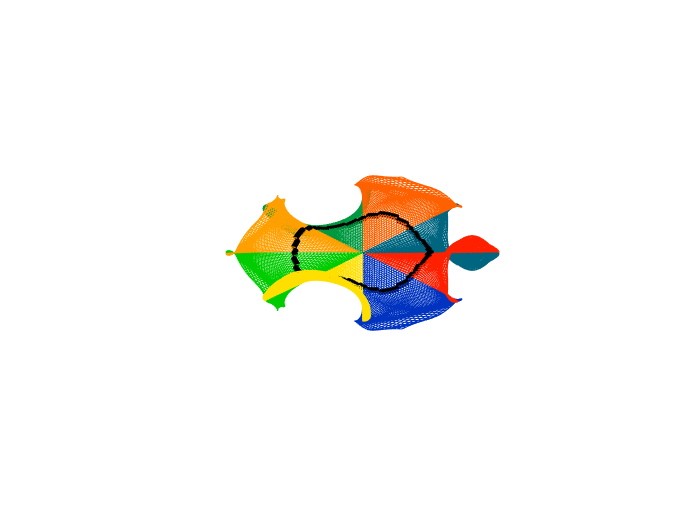}
\caption{$\eps = 1$}
\label{fig:2n4eps1}
\end{subfigure}~
\begin{subfigure}{0.31\textwidth}
\centering
\includegraphics[width=\textwidth,trim = 120 120 120 120, clip]{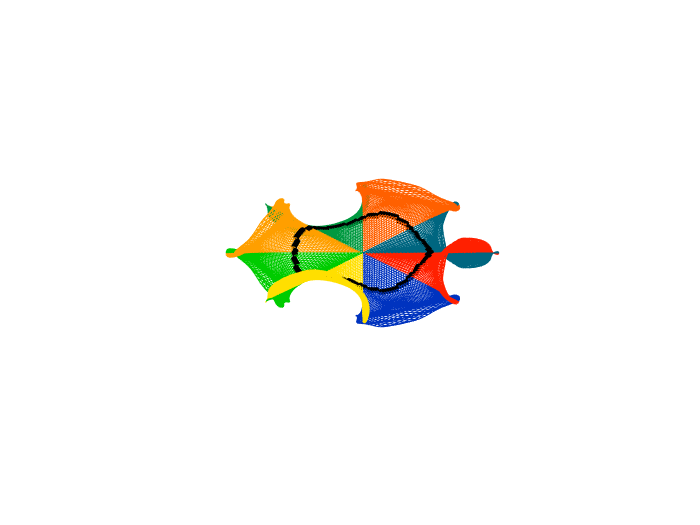}
\caption{$\eps = 2$}
\label{fig:2n4eps2}
\end{subfigure}~
\begin{subfigure}{0.31\textwidth}
\centering
\includegraphics[width=\textwidth,trim = 120 120 120 120, clip]{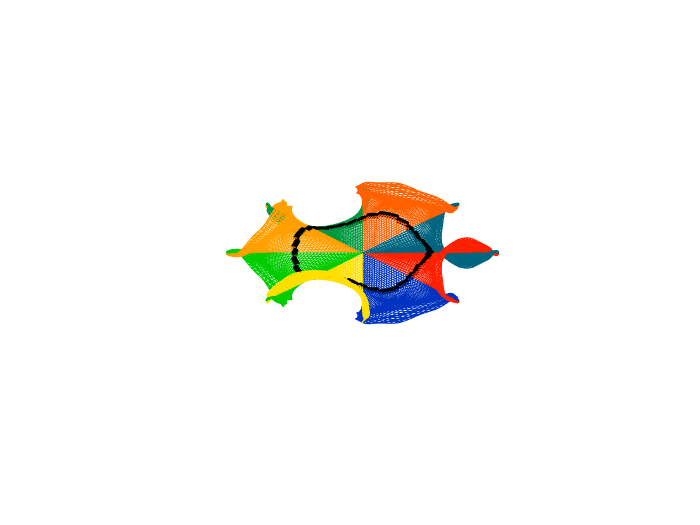}
\caption{$\eps = 3$}
\label{fig:2n4eps3}
\end{subfigure}
\caption{Amsler surfaces with 8 sectors and prescribed curvature given by \eqref{eq:PCams2}. The black ring on each surface is the set of points geodesic distance $1/2$ from the center. Note that all surfaces are the same within that ring (they are pseudospherical), and become increasingly curved outside the ring as $\eps$ increases. }
\label{fig:amslerSurfaces2}
\end{figure}

Another example is included in figure \ref{fig:amslerSurfaces2}. Here we prescribe the curvature \begin{equation}\label{eq:PCams2}K_\eps(p) = \begin{cases} -1, & D(p) \le 1/2, \\ -(1+\eps(20(D(p)-1/2))^2), & D(p) > 1/2. \end{cases} \end{equation} In this case, within the ring of geodesic distance 1/2 from the center, each of the surfaces is pseudospherical. Outside the ring of geodesic distance 1/2 from the center, the surfaces become increasingly curved as $\eps$ increases. \\

\section{Branch Points and Singular Edges} \label{sec:branch}

Note that in figure \ref{fig:amslerSurfaces1} and figure \ref{fig:amslerSurfaces2}, almost every node in each of the asymptotic complexes is ``regular," meaning that it is the point of incidence for 4 distinct quads. This is not necessarily true of the center point, which is the point of incidence for $2n$ quads when the Amsler surface is comprised of $2n$ sectors. Such irregular points are termed ``branch points." Physically, branch points allow for the type of buckling and subwrinkling which is observed, for example, in the growth of flower petals or leafy lettuce, and is not accounted for by the ``smooth" theory. These are discussed and analyzed in detail by Shearman and Venkataramani \cite{Toby}. In particular, they observe that (1) surfaces with branch points cannot result from $C^2$ isometric immersions; in fact, they cannot be approximated arbitrarily well (with respect to local bending energy) by surfaces with bounded curvature and no branch points, and (2) $C^{1,1}$ isometric immersions---which allow for branch points---can, somewhat counterintuitively, produce surfaces with significantly smaller elastic energy, despite appearing more ruffled and wrinkled to the naked eye. One reason for this is that elastic energy concentrates near singular edges, where the angle between the asymptotic lines approaches $\pi$. However, as a hyperbolic surface is evolved outward from boundary data using the Lelieuvre formulas, one can forestall the development of a singular edge by inserting branch points, and thus continue the evolution of the surface without concentrating energy. To do so, Shearman and Venkataramani \cite{Toby} develop a method of ``surgery," whereby one strategically cuts out portions of a surface and replaces the missing section with an odd number of Amsler or Pseudo-Amsler type sectors. 

Our method is also amenable to this form of surgery in order to add branch points.  To do so, one must specify which portion of the surface is to be removed, decide how exactly the removed portion is to be split into new sectors, and provide boundary data to fill in the new sectors. There are several degrees of freedom as to how these choices can be made. We describe one manner of surgery below, which maintains some symmetry regarding the newly introduced sectors, though other choices for how to build these sectors are certainly viable. 

\begin{figure}
\centering
\includegraphics[width=0.3\textwidth]{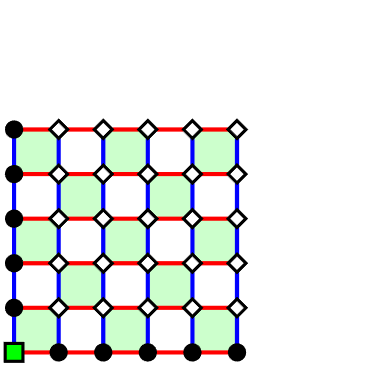} \,\,\, 
\includegraphics[width=0.3\textwidth]{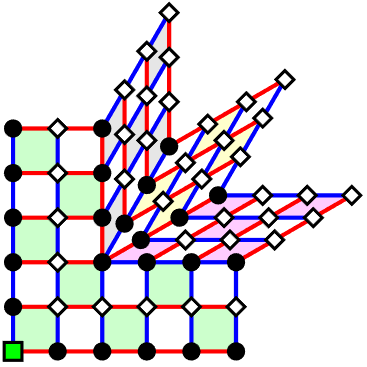} 
\caption{Beginning with an $I \times I$ asymptotic complex, we remove a square from the upper right, and replace it with $m$ copies of itself. Note that, as quad graphs, each of the parallelograms we have added are equivalent to the original square which was removed. Further note that having replaced the original square with an odd number of sectors, we can consistently label edges on the interior of the net so that each quadrilateral has two $u$-edges (red) which are opposite each other, and two $v$-edges (blue) which are opposite each other, and that the resulting quad-graph can be ``checkered" just as the original could. $s$ and $\mathbf{N}$ are prescribed on nodes marked with a Green square, the closed circles mark nodes where $\mathbf{N}$ is specified.}
\label{fig:surgNet}
\end{figure}

We suppose that one has already generated a sector of an Amsler surface in the form of a discrete immersion $\{\br_{ij}\}_{0 \le i \le I, 0 \le j \le J}$ using the method proposed in the previous section. For simplicity of exposition, we will assume that $I = J$, so that the immersion is defined on the nodes $Q = \{(i,j) \in \Z^2 \, : \, 0 \le i,j \le I\}$. We further assume that the region to be modified also a square of the form $B = \{(i,j) \in \Z^2 \, : b \le i,j \le I\}$ for some $b = 1,\ldots, I$. We ``delete" the subsurface described by $\{\br_{ij}\}_{b+1\le i,j \le I}$; this will be replaced with collection of new Amsler sectors given by immersions $\br^{(k)}: B \to \R^3$ and corresponding normal fields $\bN^{(k)}:B \to \mathbb S^2$, for $k = 1,\ldots, m$. In doing this, we are replacing the section $B$ of the asymptotic complex $Q$ with $m$ ``copies" of itself. This is demonstrated with a $5\times 5$ asymptotic complex in figure \ref{fig:surgNet}. In the figure, we use $b = 2$ and $m = 3$, so that a $3 \times 3$ square is removed and replaced with $3$ copies of itself. It is proven in \cite{Toby}, that in order for the resulting quad graph to still be an asymptotic complex, $m$ must be odd. Otherwise, there will be no consistent way to label edges so that each quadrilateral consists of two $u$-edges which are opposite each other, and two $v$-edges which are opposite each other. In figure \ref{fig:surgNet}, the pink-tiled parallelogram is the domain for $\br^{(1)}$ which shares a $u$-boundary with the original complex, the yellow-tiled parallelogram is the domain for $\br^{(2)}$, and the gray-tiled parallelogram is the domain for $\br^{(3)}$ which shares a $v$-boundary with the original parallelogram. We emphasize here that while the new parallelograms are equivalent as graphs to the excised square, the orientation of alternate parallelograms is flipped, in that the ``right-hand" boundary is alternately comprised of $u$-edges or $v$-edges. 

\begin{figure}[b!]
\centering
\includegraphics[width=0.3\textwidth,trim=150 50 150 200,clip]{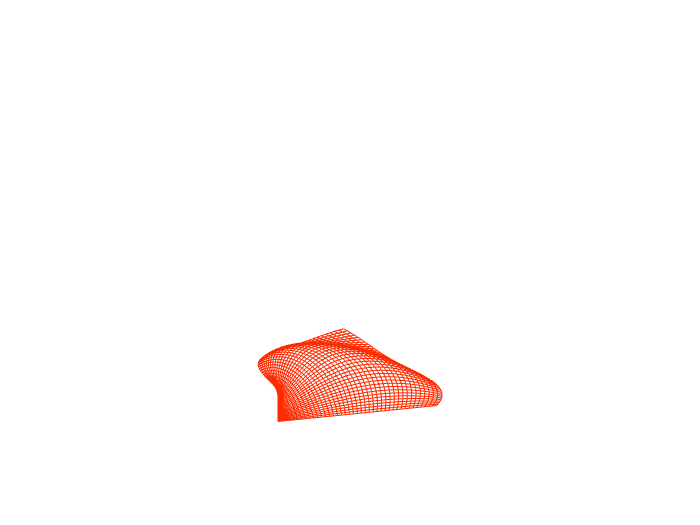} 
\includegraphics[width=0.3\textwidth,trim=150 50 150 200,clip]{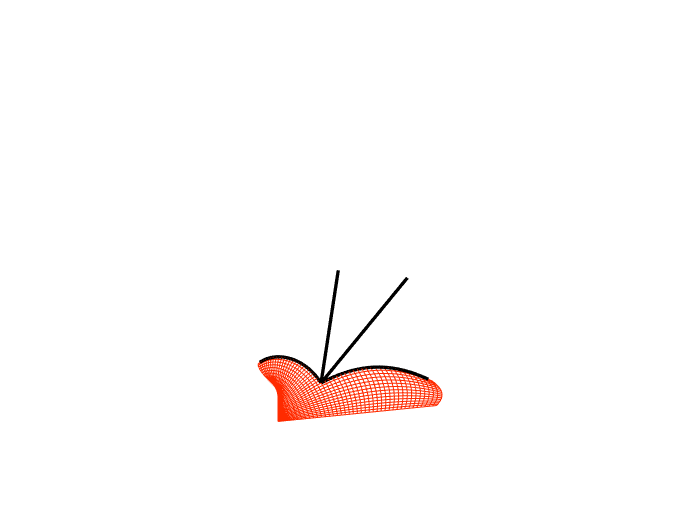} 
\includegraphics[width=0.3\textwidth,trim=150 50 150 200,clip]{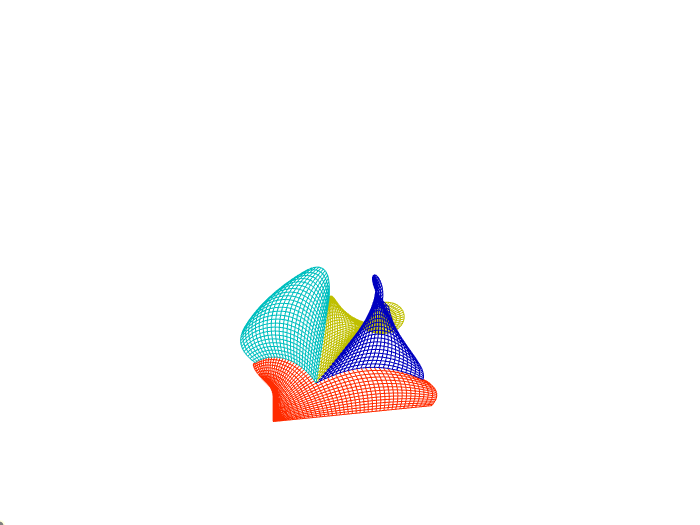} 
\caption{A visualization of ``surgery" to introduce a branch point. Beginning with an Amsler sector, we cut out a square, split the angle into $m$ equal parts, specify boundary data along the new axes (black), and then generate new Amsler sectors using the specified boundary data along the black lines, and the data from the original sector along the boundaries of the cut. The two sectors that are bounded by a straight line and a space curve from the initial immersion are Pseudo-Amsler sectors, while the ``middle" sector bounded by two straight lines is an Amsler sector. }
\label{fig:branchPoint}
\end{figure}

In general, we define these new sectors so that $\br^{(1)}_{ib} = \br_{ib}$ for all $b \le i \le I$ and $\br^{(m)}_{bj} = \br_{bj}$ for all $b \le j \le I$, and similarly for the normal fields $\bN^{(1)}$ and $\bN^{(m)}$. That is, $\br^{(1)}$ uses the values of the original immersion as $u$-boundary data, and $\br^{(m)}$ uses the values of the original immersion as $v$-boundary data. The other sectors are counted counterclockwise with the $\br^{k}_{bi} = \br^{k+1}_{ib}$ for $b \le i \le I$ and $k = 1,\ldots m-1$. It remains to specify these boundary values, and the corresponding normal field values. To do so, we consider the behavior of the original immersion near the branch point. Specifically, we consider the plane formed by $\br_{b+1,b} - \br_{bb}$ and $\br_{b,b+1} - \br_{bb}$, and let $\theta_b$ be the angle between these vectors. In our construction, the boundaries of the sectors $\br^{(k)}$ are the straight lines in this plane which split $\theta_b$ into $m$ equal parts. Since the normal to this plane is $\bN_{bb}$, we take the ``left-hand" boundary of $\br^{(k)}$ for $k=1,\ldots,m-1$ to be a straight line proceeding from $\br_{bb}$ in the direction of $$\text{Rot}\left(\frac{k\theta_b}m; \bN_{bb}\right)(\br_{b+1,b} - \br_{bb}).$$ One can then specify units of distance to initialize $\br^{(k)}$ along these lines, and initialize the normal fields along these lines exactly as described above for the initial Amsler surface. Having done so, one can solve the discrete non-constant curvature Lelieuvre equations on each sector exactly as described above. We demonstrate this procedure in figure \ref{fig:branchPoint}. Note that the original sector (figure \ref{fig:branchPoint}, left) has a singular edge where the surface curls over on itself. We excise a square containing this edge, specify new axes and new boundary data in the missing sector, (figure \ref{fig:branchPoint}, center), and fill in the missing sector with a collection of new (Pseudo-)Amsler sectors (figure \ref{fig:branchPoint}, right).

\begin{figure}[t!]
\centering
\includegraphics[width=0.48\textwidth,trim=150 50 150 200,clip]{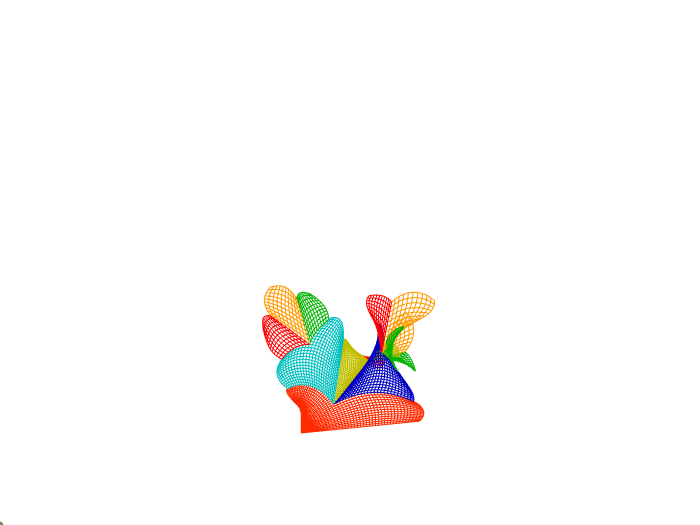}
\includegraphics[width=0.48\textwidth,trim=125 50 150 100,clip]{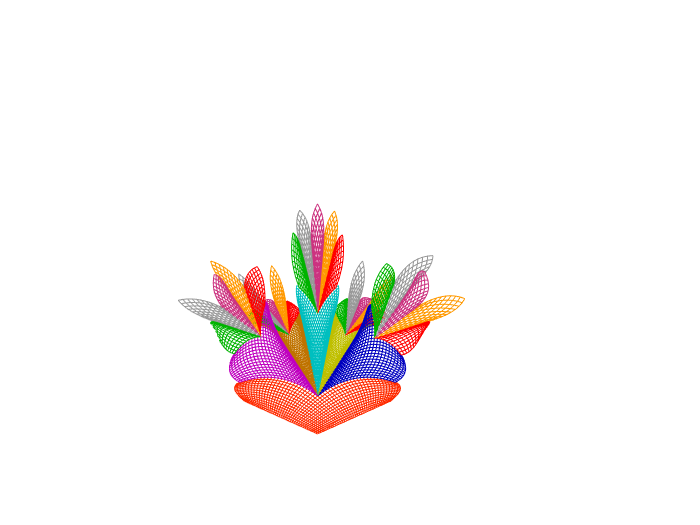}
\caption{Left: The discrete hyperbolic surface from figure \ref{fig:branchPoint} with additional branch points introduced in the new sectors. Right: The same discrete surface with 5 new sectors introduced at each branch point instead of 3.}
\label{fig:branchPointTwoGen}
\end{figure}

This procedure can then be recursively applied to insert new branch points by excising pieces of the new (Pseudo-)Amsler sectors. The astute reader will note that the new sectors introduced through the branch point in figure \ref{fig:branchPoint} still have a singular edge (or at least a very nearly singular edge). Inserting another branch point in each of the new sectors will ameliorate this. This is demonstrated in figure \ref{fig:branchPointTwoGen}, where we introduce another generation of branch points, and the newly introduced sectors do not have singular edges. 

\begin{figure}[t!]
\begin{subfigure}{0.3\textwidth}
\centering 
\includegraphics[width=\textwidth,trim = 80 40 80 40, clip]{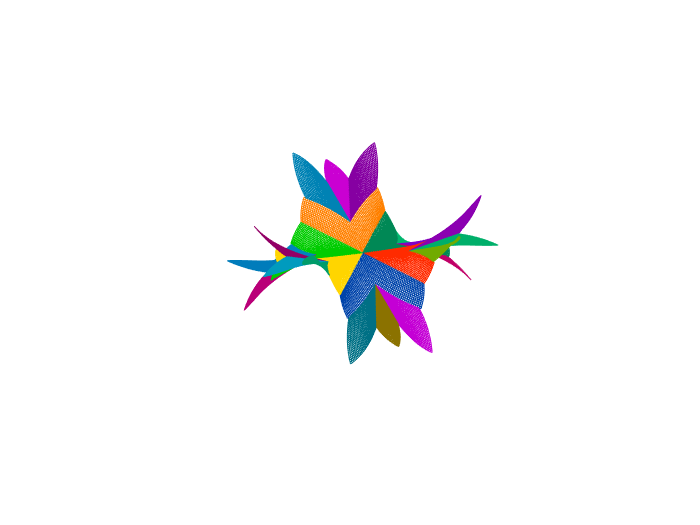}
\caption{$m=3, \eps = 1$}
\label{fig:n3m3eps1}
\end{subfigure}~
\begin{subfigure}{0.3\textwidth}
\centering 
\includegraphics[width=\textwidth,trim = 80 40 80 40, clip]{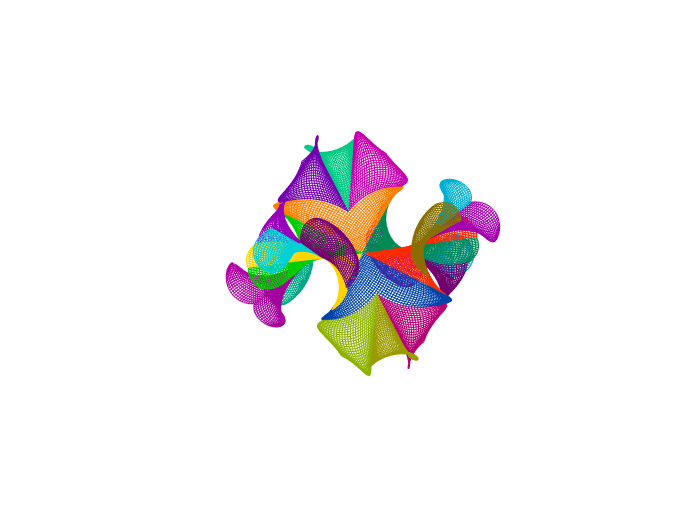}
\caption{$m=3, \eps = 10$}
\label{fig:n3m3eps10}
\end{subfigure}~
\begin{subfigure}{0.3\textwidth}
\centering 
\includegraphics[width=\textwidth,trim = 80 40 80 40, clip]{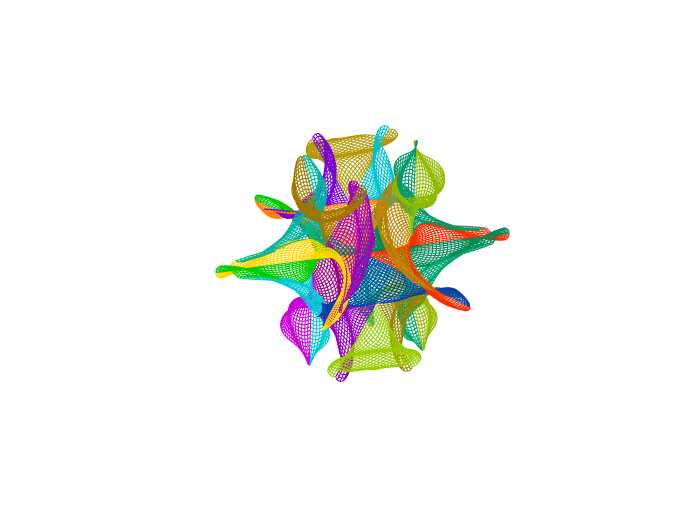}
\caption{$m=3, \eps = 50$}
\label{fig:n3m3eps50}
\end{subfigure}~\\
\begin{subfigure}{0.3\textwidth}
\centering 
\includegraphics[width=\textwidth,trim = 80 40 80 40, clip]{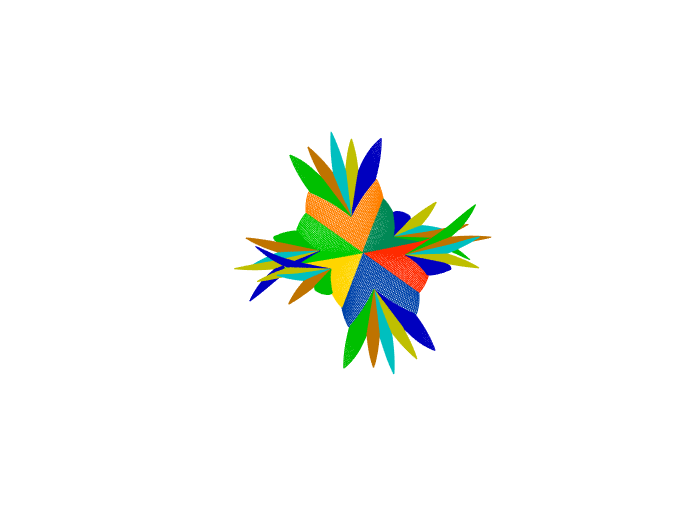}
\caption{$m=5, \eps = 1$}
\label{fig:n3m5eps1}
\end{subfigure}~
\begin{subfigure}{0.3\textwidth}
\centering 
\includegraphics[width=\textwidth,trim = 80 40 80 40, clip]{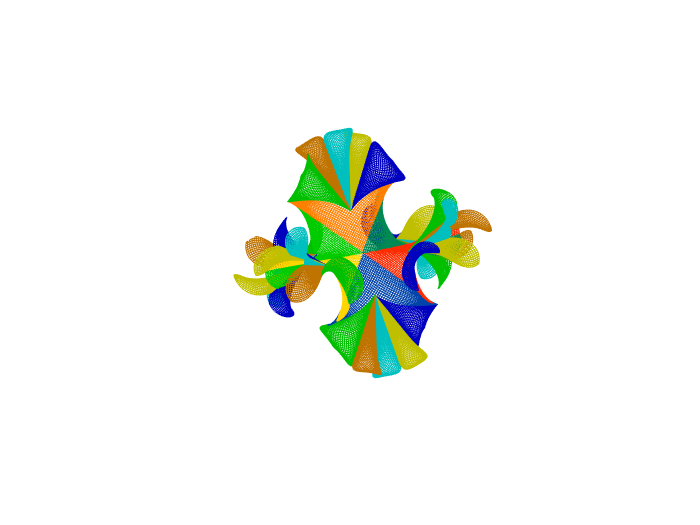}
\caption{$m=5, \eps = 10$}
\label{fig:n3m5eps10}
\end{subfigure}~
\begin{subfigure}{0.3\textwidth}
\centering 
\includegraphics[width=\textwidth,trim = 80 40 80 40, clip]{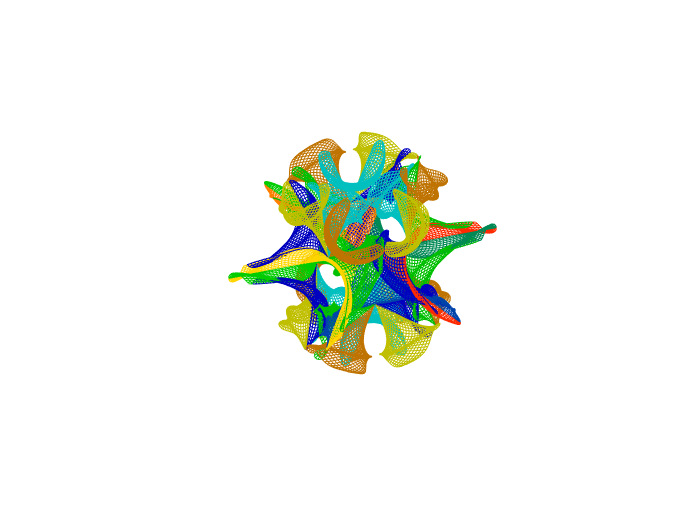}
\caption{$m=5, \eps = 50$}
\label{fig:n3m5eps50}
\end{subfigure}~
\caption{The Amsler surfaces from figure \ref{fig:amslerSurfaces1}(d)-(f) with one generation of branch points introduced. Here the prescribed curvature is $K(p) = - (1+\eps D(p))$, where $D(p)$ is geodesic distance to the origin. Each surface is comprised of $6$ Amsler sectors and we introduce $m = 3$ or $m = 5$ new sectors at the branch point. }
\label{fig:AmslerSurfacesFullBP}
\end{figure}

While the demonstrations in figure \ref{fig:branchPoint} and figure \ref{fig:branchPointTwoGen} use only a single sector, these can pasted together as before to create closed surfaces whose normal field is continuous across the boundary lines. In figure \ref{fig:AmslerSurfacesFullBP}, we generate the same discrete surfaces as in figure \ref{fig:amslerSurfaces1}(d)-(f), but with branch points added in each sector. At the branch points, we alternately introduce $m = 3$ or $m=5$ new sectors. 

\section{The Strip Geometry: A Model for Edges}
\label{sec:strip_geometry}

In the previous section, we considered a disk-like geometry where curvature depends on the geodesic distance from a central point, a model well-suited for objects that grow isotropically. However, many of the phenomena that motivate our work, such as the frilly edges of leaves or the crenellations on torn plastic sheets \cite{Sharon2004Leaves, Plastic1}, are better described by an asymptotically flat strip with a decorated edge. To model this, we construct an immersion using a composite curvature function that combines a power-law behavior near the edge with an asymptotically flat far-field.

\subsection{Numerical Construction of the Strip Immersion}

Our numerical construction of the strip immersion follows a three-step procedure. First, we define an asymptotic coordinate system using an analytically tractable {\em seed surface}. Second, we define our target {\em intrinsic curvature model}, which is a composite function designed to replicate experimental observations. Finally, we use the seed surface to provide boundary data for an {\em iterative algorithm} that deforms the initial surface until its geometry matches the target curvature.

\begin{figure}[t!]
    \centering
    \includegraphics[width=0.32\textwidth]{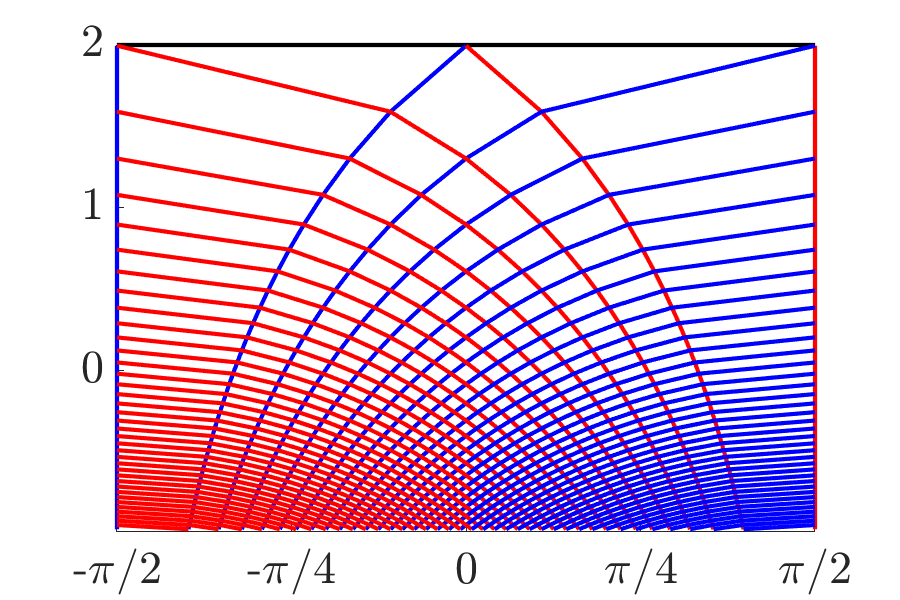}\, 
    \includegraphics[width=0.4\textwidth]{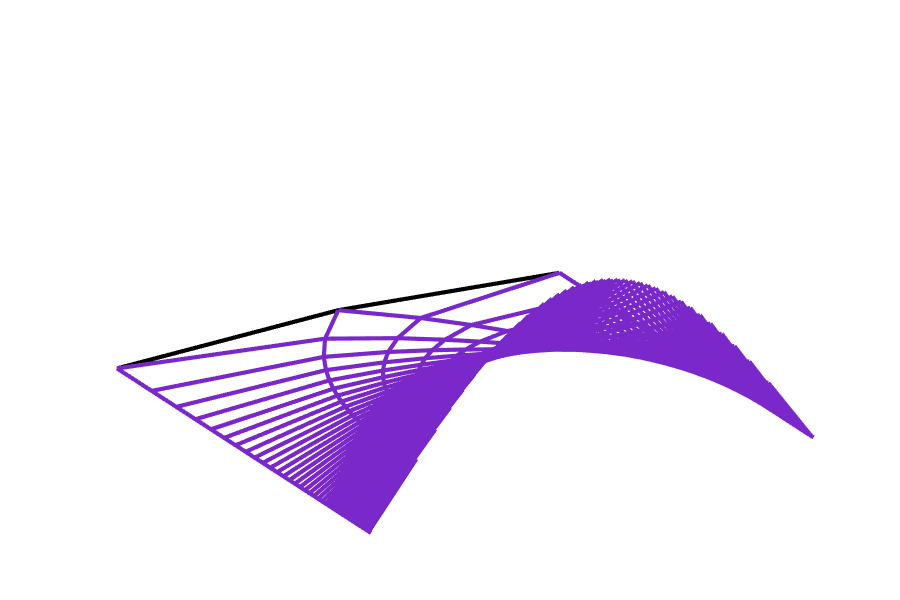}\,
    \includegraphics[width=0.24\textwidth]{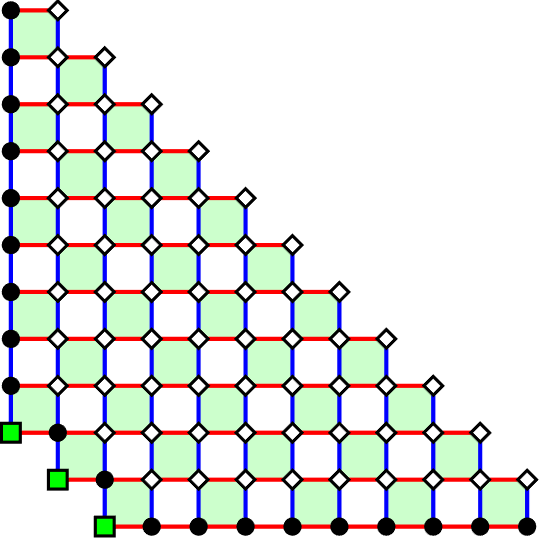}
    \caption{Left: 2D projection of asymptotic coordinates on the seed surface in $\Omega = [-\pi/2,\pi/2] \times [0,2]$. Center: The seed surface plotted in these coordinates. Right: The quadgraph of the computational domain. The bottom-left diagonal corresponds to the far boundary ($y=2$). The figure illustrates the Goursat problem setup: data for the normal $\mathbf{N}$ is specified on the boundary (black circles), while both $\mathbf{N}$ and the geodesic distance $s$ are specified at the corner (green squares).}
    \label{fig:strip_immersion}
\end{figure}

We begin with a seed surface given by a Monge patch $(x,y,w(x,y))$ where $(x,y) \in \Omega \defeq (-\tfrac \pi 2, \tfrac \pi 2) \times [0, 2]$ and $w(x,y) = Ae^{-y}\cos(x)$ with the constant $A$ to be specified later. This seed surface has Gauss curvature given by
\begin{equation}\label{eq:curvStripExplicit} 
K = \frac{w_{xx}w_{yy} - w_{xy}^2}{(1 + w_x^2 + w_y^2)^2} = \frac{-A^2 e^{-2y}}{(1 + A^2 e^{-2y})^2}. 
\end{equation} 
In particular, the curvature depends only on $y$, is always negative $K(y) < 0$, and approaches 0 as $y$ increases, so the surface is nearly flat near the ``upper" boundary at $y = 2$. Since $K <0$, we can foliate the surface by asymptotic curves that are given by \begin{equation} \label{eq:asympCurvesEqnStrip}
    -\cos(x)\,dx^2 + 2\sin(x)\,dxdy + \cos(x)\,dy^2 = 0.
\end{equation} 
Integrating gives two families of curves, from which we define the asymptotic coordinates $(u,v)$:
\begin{equation}\label{eq:asympCoordsStrip}
\begin{aligned} 
u &= e^{-y}(1+\sin(x)),\\ 
v &= e^{-y}(1-\sin(x)), 
\end{aligned} 
\quad\quad \leftrightsquigarrow \quad\quad 
\begin{aligned} 
x &= \arcsin\left(\tfrac{u-v}{u+v}\right), \\ 
y &= -\log\left(\tfrac{u+v}{2}\right), \\ 
z & = A \sqrt{uv}. 
\end{aligned}
\end{equation} 
The left and right boundaries $x=\pm \pi/2$ are themselves asymptotic curves where $u=0$ and $v=0$, respectively. A visualization of these coordinates on the seed surface is shown in Figure~\ref{fig:strip_immersion} (left and center panels).

This coordinate system suggests a computational domain equivalent to a truncated square grid, shown in Figure~\ref{fig:strip_immersion} (right panel). The discrete vertices $(i,j)$ with $i+j = M_{\text{min}}$ correspond to the line $y = 2$, and those with $i+j=M_{\text{max}}$ correspond to $y= 0$. This setup is ideal for a Goursat-type problem where initial data is specified on the ``past" boundaries (black and green points) and the solution is evolved forward.

While the seed surface provides a useful coordinate system, its curvature depends on the \emph{extrinsic} coordinate $y$. To create a model based on the material's internal geometry, we instead define an \emph{intrinsic} curvature that depends on $s$, the geodesic distance from the wrinkled edge. We prescribe the total geodesic width of the strip to be 2, so the known far boundary corresponds to $s=2$.

To replicate the geometry of torn plastic sheets, we construct a composite function, $-K_{\text{comp}}(s)$, that combines two distinct behaviors. Near the wrinkled edge, we use a power-law form, $-K_2(s) \propto (s+\ell)^{-2}$, consistent with experimental observations \cite{Sharon2004Leaves,Plastic1}. For the far-field, we adopt the asymptotically flat form from the seed surface's curvature, $-K_1(s)$, by replacing the role of extrinsic coordinate $y$ with the intrinsic distance $s$. By choosing parameters $A = e/\sqrt{2}$, $\ell = 2$, and a contact point of $s_0 = 1$, these two forms are joined into the following $C^1$ function:
\begin{equation}
    \label{eq:compositeCurve}
    -K_{\text{comp}}(s) = \begin{cases} 
    \frac{2}{(s+2)^2}, & s < 1 \\
    \frac{2\exp(2 - 2s)}{(2+\exp(2-2s))^2}, & s \ge 1
    \end{cases}
\end{equation}

\begin{figure}[t!]
    \centering
    \includegraphics[width=0.4\textwidth]{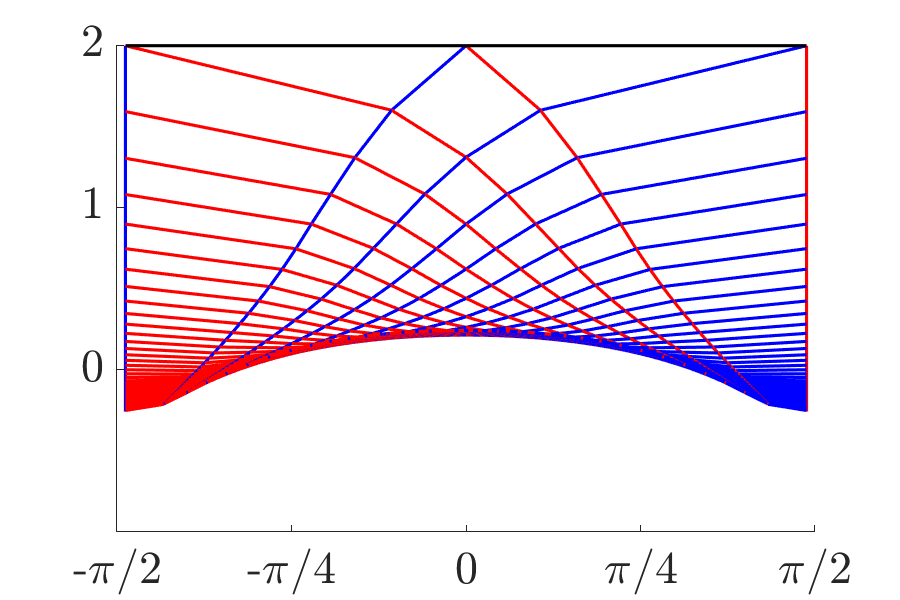}\,
    \includegraphics[width=0.5\textwidth,trim = 0 10 0 60,clip]{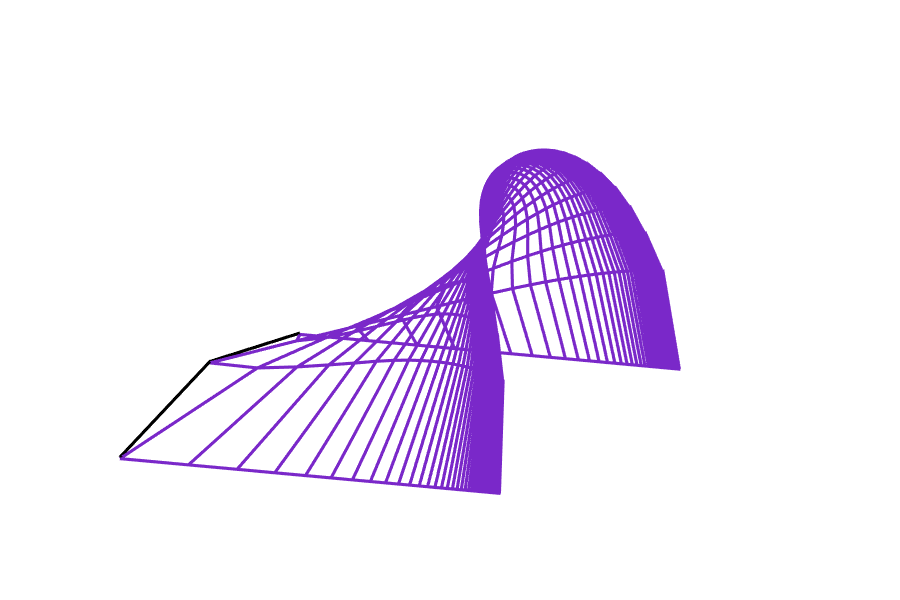}
    \caption{Two views of the discrete surface with the intrinsic curvature $K_{\text{comp}}(s)$ given by \eqref{eq:compositeCurve}. There is a singular edge where the $u$- and $v$-asymptotic curves touch tangentially and the surface is no longer the graph of a differentiable function.}
    \label{fig:stripOrig}
\end{figure}

To find the immersion for this target curvature, we use the seed surface to provide boundary data and an initial guess, $\B r^{(0)}_{ij}, \B N^{(0)}_{ij}$. In particular, the seed surface, with curvature $K_1(y)$ is very close to an immersion of surface with curvature $K_1(s)$ as the approximation $y \approx s$ is consistent in the far field. We then run Algorithm~\ref{alg:IterationAlg}, with the modification that the fast marching method decreases $s$ starting with $s=2$ on the edge $y=2$, to update the immersion based on $K_{\text{comp}}(s)$, until the process converges. The result of this procedure, shown in Figure~\ref{fig:stripOrig}, is a surface that develops a singular cuspidal edge where the asymptotic curves from the $u$- and $v$- families touch tangentially.

\begin{figure}[t!]
    \centering
    \hspace{0.3cm}\includegraphics[width=0.25\textwidth,angle = 225,origin=c]{stripGrid.png}\hspace{1cm} 
    \includegraphics[width=0.25\textwidth,angle = 225,origin=c]{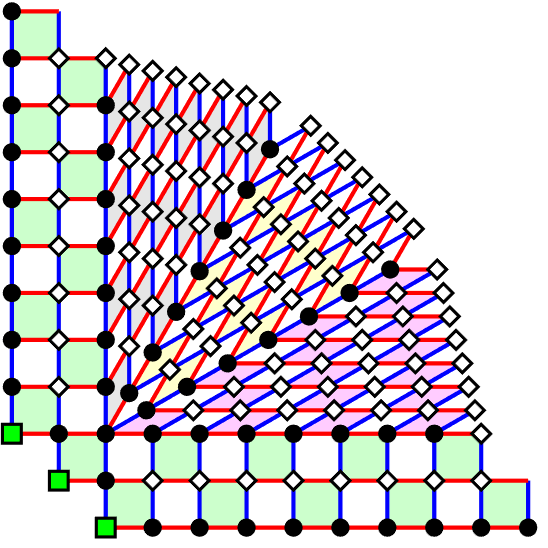}
    \vspace{-2cm}
    \caption{Introducing a branch point via surgery removes the singular edge. The grid on the left represents the computational domain for the unbranched surface depicted in Fig.~\protect{\ref{fig:stripOrig}} that has a singular edge. 
    On the right is a visualization of the surgery on the discrete computational domain, shown here for the case of 3-fold subwrinkling surface depicted in Fig.~\protect{\ref{fig:strip_surgery2}}.}
    \label{fig:strip_surgery}
\end{figure}

\begin{figure}[t!]
    \centering
    \includegraphics[width=0.45\textwidth]{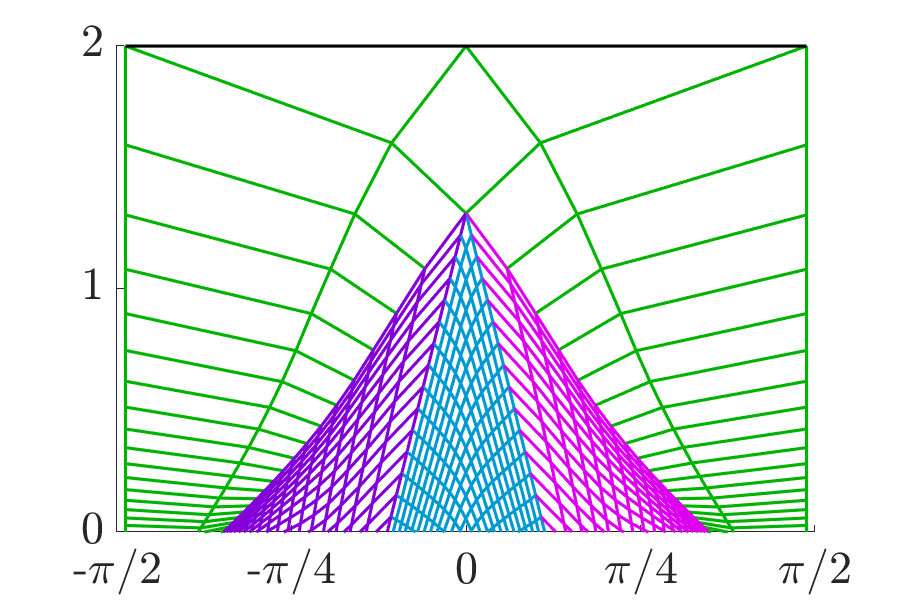} \includegraphics[width=0.45\textwidth]{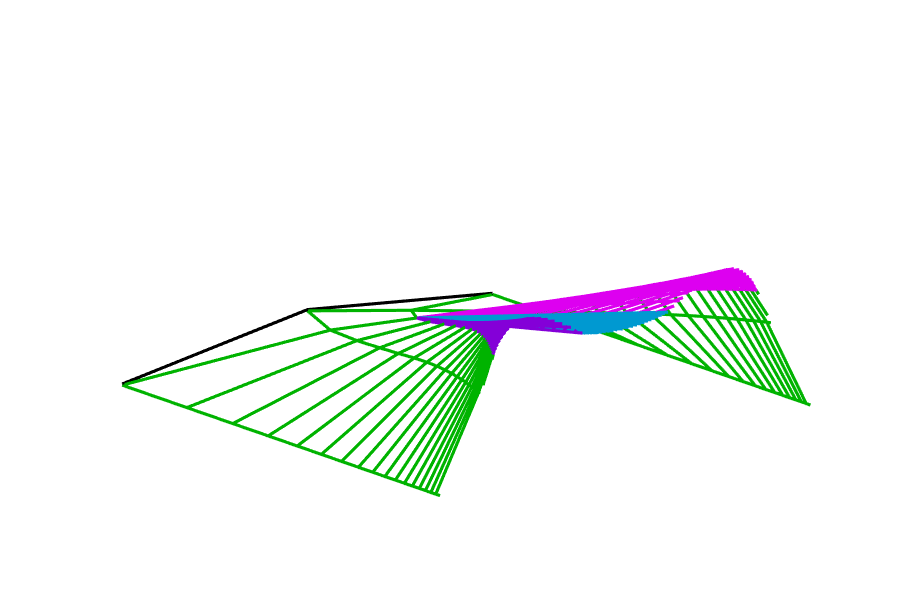} \\
    \includegraphics[width=0.45\textwidth]{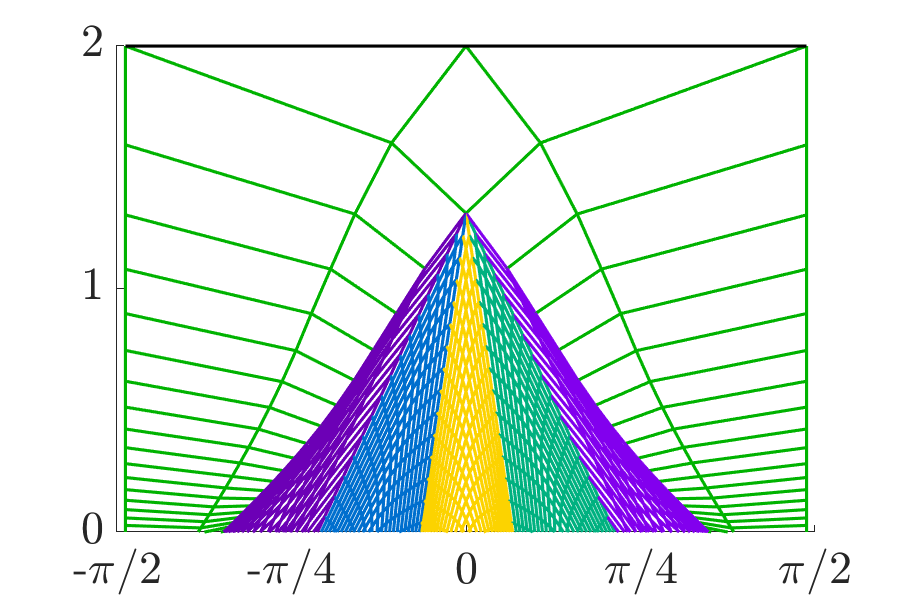} \includegraphics[width=0.45\textwidth]{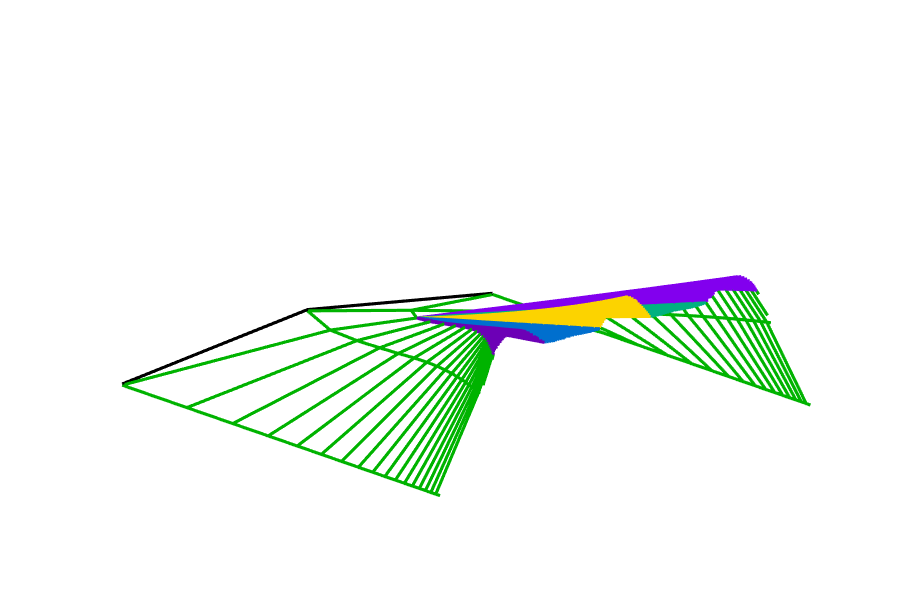} \\
    \includegraphics[width=0.45\textwidth]{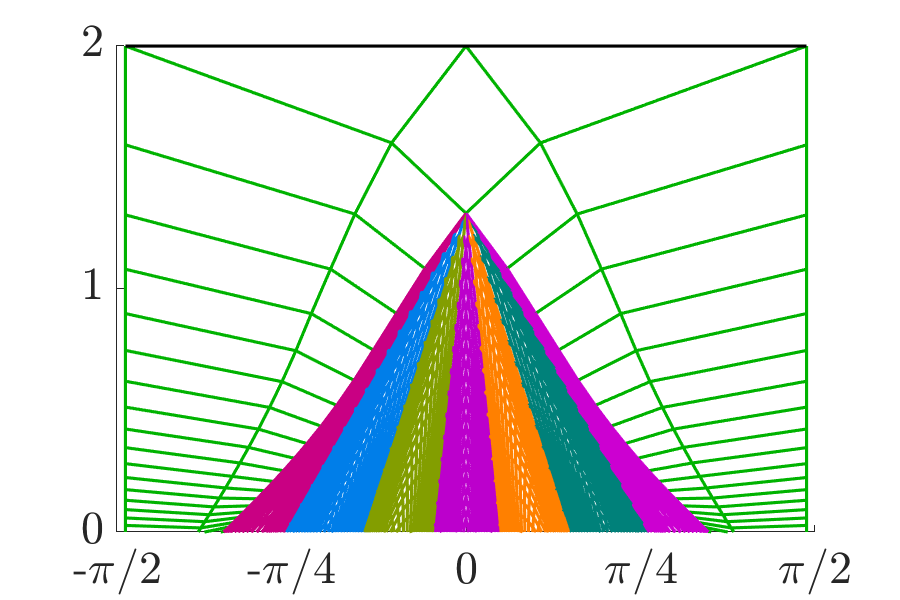}  \includegraphics[width=0.45\textwidth]{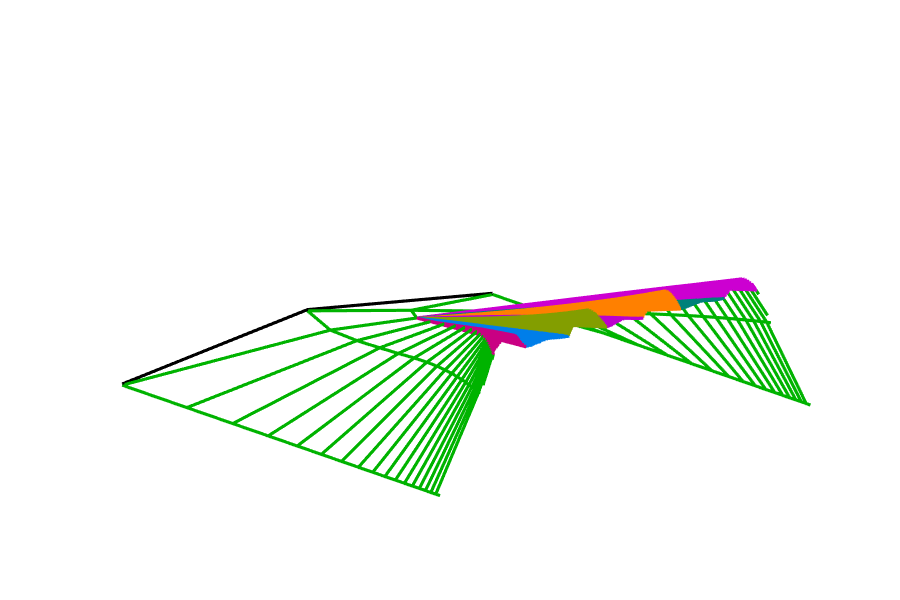} 
    \caption{Surfaces with curvature given by $K_{\text{comp}}(s)$, where $s$ is the geodesic distance to the wrinkled edge. The far boundary is normalized to $s=2$. We have introduced branch points correponding to 3-fold (top), 5-fold (middle), and 7-fold (bottom) subwrinkling.}
    \label{fig:strip_surgery2}
\end{figure}

As with the Amsler-type surfaces, the bending energy concentrates along this singular edge, which is not physically preferred. We therefore perform surgery to introduce branch points and forestall this singularity. The procedure is analogous to that of the previous section: we choose a cut point $b$, remove grid nodes $(i,j)$ with $i,j\ge b$, and replace the excised portion with $m$ copies, pasted along new boundary lines (Figure~\ref{fig:strip_surgery}). The final branched surfaces, shown in Figure~\ref{fig:strip_surgery2} for 
3-fold, 5-fold and 7-fold subwrinkling have no singular edge. The surgically inserted regions successfully demonstrate the sort of high-frequency wrinkling observed at the edge of a torn plastic strip \cite{Sharon2004Leaves,Plastic1}.

\section{Conclusion}

 Surfaces of non-constant negative curvature can serve as mathematical models for thin elastic sheets that arise in countless physical scenarios, both natural and artificial. In this paper, we describe a method for generating surfaces of prescribed non-constant negative curvature {\em with distributed branch points}. 
 The key idea is to use the non-constant curvature version of the classical Lelieuvre formulas in both a continuous and discrete setting. In the case of non-constant curvature, the Lelieuvre formulas are no longer fully explicit, so we devise an iterative algorithm to numerically integrate them. We are especially interested in scenarios where curvature is a function of ``material" coordinates, in particular the geodesic distance to some point or edge. Accordingly, we describe a fast marching method for computing geodesic distance on a triangulated manifold. Using our method, we generate branched hyperbolic surfaces, which are comprised of bounded sectors patched together so as to maintain a continuous normal map throughout the surface. Finally, we demonstrate that our method is amenable to a sort of `surgery' to introduce branch points, whereupon sectors of a surface are excised and replaced by several new sectors. This allows for surfaces of non-constant negative curvature which exhibit the type of multi-generational buckling and wrinkling often observed in physical applications. 

 There is a larger context for this work. Motivated by applications, there is interest in developing numerical methods for solving coupled systems of PDEs on (possibly evolving) surfaces. The typical approach has been through finite elements defined on a Lagrangian triangulation, i.e. a mesh that is fixed in the material, since physical processes on the surface are often naturally modeled through PDEs in material coordinates \cite{Bartels2012FiniteEM,dziuk2013finite}. While successful, this approach does not exploit the particular simplifications that arise by having a grid that is adapted to the geometry, i.e. one that discretizes the asymptotic curves on the surface. 
 
\begin{figure}
     \centering
     \includegraphics[width=0.6\textwidth,trim=150 40 50 150, clip]{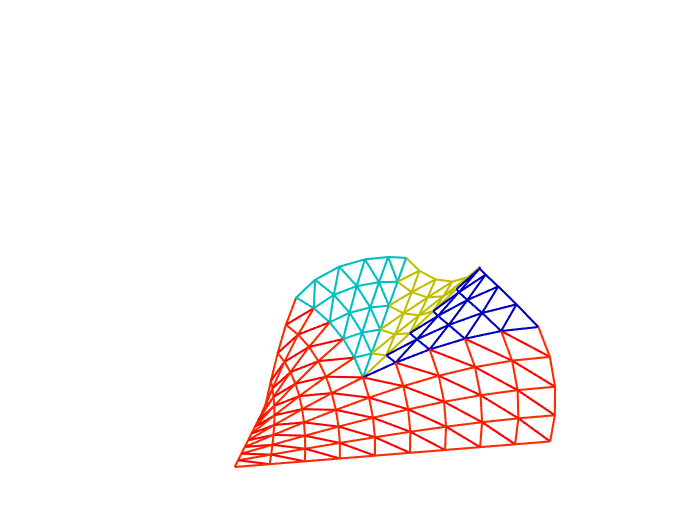}
     \caption{A discretization of surface like that in figure \ref{fig:branchPoint} which respects both the ``geometry" (by discretizing asymptotic curves) and the ``physics" (by including a triangulation with no obtuse angles).}
     \label{fig:FEM}
 \end{figure}

 We have considered a prototypical problem of this type, coupling ``geometry", namely the prescribed Gauss curvature equation, with ``physics", in this case the eikonal equation that determines the geodesic distance. To solve such coupled systems numerically, especially in scenarios involving potential branch points, we expect that a discretization that respects the extrinsic geometry of the surface and is ``branch-point aware" will be more efficient. The underlying asymptotic complex for a hyperbolic surface does provide such a discretization and the mapping of this graph into the surface $\Omega$ gives a cellular decomposition into hyperbolic quadrilaterals. In particular, the sum of the angles in each quadrilateral is less than $2 \pi$. This allows us to define a geometrically ``natural" triangulation of $\Omega$ by bisecting quadrilaterals into 2 triangles, 
 neither of which has an obtuse angle, thus generating a good-quality mesh. This procedure of obtaining a natural triangulation of $\Omega$ is illustrated in figure \ref{fig:FEM}. We expect that finite elements based on this approach will play an important role in solving coupled surface equations on evolving hyperbolic surfaces, especially in the presence of branch points.

\section*{Acknowledgments}

The authors were supported in part by the Research Training Group in Data Driven Discovery at the University of Arizona (NSF award DMS-1937229) and by NSF award DMR-1923922. SV was supported in part by NSF awards DMS-2108124 and DMS- 2511503.

 \appendix 
 \section{Proof of Convergence of Algorithm \ref{alg:IterationAlg} for Amsler-Type Surfaces} \label{sec:convProof} In this appendix, we prove the convergence of Algorithm \ref{alg:IterationAlg} which generates a discrete surface of Gaussian curvature $K_\eps(p)$ which is seen as an size-$\eps$ perturbation of given discrete surface with curvature $K(p)$. The overarching goal is to prove that the steps in Algorithm \ref{alg:IterationAlg} which lead from $(\br^{(n)}_{ij}, \bN^{(n)}_{ij})$ to $(\br^{(n+1)}_{ij}, \bN^{(n+1)}_{ij})$, when taken holistically, are contractive given that $\eps$ is small enough, and thus the iteration converges via the Banach Fixed Point Theorem. Recall, given $(\br^{(n)}_{ij},\bN^{(n)}_{ij})$, the steps constituting a single iteration of Algorithm \ref{alg:IterationAlg} are as follows: \begin{itemize}
     \item[(I)] Compute the discrete geodesic distance $D^{(n)}_{ij}$ function which gives distance from $\br^{(n)}_{ij}$ to the origin.
     \item[(II)] Set $K^{(n)}_{ij} = K_{ij} - \eps h(D^{(n)}_{ij})$
     \item[(III)] Compute $(\br^{(n+1)}_{ij}, \bN^{(n+1)}_{ij})$ using \eqref{eq:fullUpdate} with curvature $K^{(n)}_{ij}$ (after having initialized the boundaries of $\br^{(n+1)}_{ij}$ using the same boundary data as $\br^{(0)}_{ij}$)
 \end{itemize} 

We prove error bounds for steps (I) and (III) in lemmas before we proceed with the proof of convergence. 

Before we commence with the proof, we note that, while this proof only directly applies to the class of surfaces we are interested in (the Amsler-type surfaces described in section \ref{sec:amsler}), with small modifications, it could likely be extended to a much broader class of surfaces. The key piece of information to kickstart the algorithm is the manner in which boundary data is specified so that the discrete Leliuevre equations can be resolved. At minimum, boundary data should be specified along discrete sets $\Gamma$ with the following property: for each vertex $V \not\in \Gamma$, there is precisely one discrete $u$-asymptotic curve and one discrete $v$-asymptotic curve each of which intersect $V$ and $\Gamma$. This property would allow one to march from any given vertex to the boundary and thus accomplish bounds like those in \eqref{eq:NijBoundFinal} and \eqref{eq:rijBound}. Along with this, one may like to require that $\Gamma$ is as small as possible in some sense, upon which $\Gamma$ could be dubbed an \emph{asymptotic skeleton}. Investigation into different types of such sets, as well as the topological differences in the resulting surfaces, is reserved for future work. 

\begin{lemma} \label{lem:1}
Assume that $\{\br_{ij}\}$ and $\{\br_{ij}^*\}$ for $i = 0,1,\ldots, I$ and $j=0,1\ldots,J$  are generated from the Lelieuvre equations \eqref{eq:fullUpdate} and that $\{D_{ij}\}$ and $\{D^*_{ij}\}$ are discrete functions measuring geodesic distance for $\{\br_{ij}\}$ and $\{\br_{ij}^*\}$ respectively. Then there is a constant $C>0$ such that $$\max_{1 \le i \le I, 1 \le j \le J} \abs{D_{ij}-D^*_{ij}} \le C \left(\max_{1 \le i \le I, 1 \le j \le J}  \|\br_{ij} - \br_{ij}^*\| \right).$$  \\
\end{lemma}

 \emph{Proof.} From the discrete manifold $\{\br_{ij}\}$,  we define a manifold $\br$ on a continuum by zooming into a quad with vertices $\br_{ij}, \br_{i+1,j}, \br_{i,j+1}, \br_{i+1,j+1}$ and setting $$\br(u,v) = (i+1-u)(j+1-v)\br_{ij} + (i+1-u)(v-j)\br_{i+1,j} + (u-i)(j+1-v)\br_{i,j+1} +  (u-i)(v-j) \br_{i+1,j+1},$$ where $(u,v) \in [i,i+1]\times [j,j+1]$ are local coordinates on the quad for $i=0,1,\ldots, I-1, j = 0,1,\ldots, J-1$. We likewise define $\br^*$ from $\br_{ij}^*$. Having done so, we have two different manifolds which are defined on the same underlying coordinates $(u,v) \in [0,I] \times [0,J]$. The manifolds have respective Riemannian metrics  \begin{equation} \begin{split} 
g &= \|\br_u\|^2du^2 + 2\langle \br_u,\br_v\rangle dudv + \|\br_v\|^2 dv^2, \\
g^* &= \|\br^*_u\|^2du^2 + 2\langle \br^*_u,\br^*_v\rangle dudv + \|\br^*_v\|^2 dv^2.
\end{split} \end{equation} These derivatives are expressed entirely in terms of the discrete coordinates: \begin{equation} \label{eq:rderiv}\begin{split} \br_u &= (j+1-v)(\br_{i,j+1}-\br_{ij}) + (v-j)(\br_{i+1,j+1} - \br_{i+1,j}), \\ 
\br_v &= (i+1-u)(\br_{i+1,j} - \br_{ij}) + (u-i)(\br_{i+1,j+1} - \br_{i,j+1}), \end{split}  \end{equation}  when $(u,v) \in [i,i+1] \times [j,j+1]$, and similarly for $\br^*$.  From this, using that $$ \br_{i+1,j} - \br_{ij} = \bnu_{i+1,j}\times \bnu_{ij} \,\,\,\,\, \text{ and } \,\,\,\,\, \br_{i,j+1} - \br_{ij} = - \bnu_{i,j+1}\times \bnu_{ij}$$ as well as $\|\bnu_{ij}\| \le 1$ (which holds as long as $\br_{ij}$ was generated using a curvature function satisfying $\abs{K(p)}\ge 1$), we immediately have $\|\br_u\|, \|\br_v\| \le 1,$ and likewise for $\br^*$.  Then \begin{equation}  \begin{split}
\abs{\|\br_u\|^2 - \|\br^*_u\|^2} &= \abs{\innerprod{\br_u + \br_u^*}{\br_u -\br_u^*}} \le 2 \|\br_u - \br_u^*\|. \end{split}
\end{equation} From \eqref{eq:rderiv}, we see that $$\|\br_u - \br_u^*\| \le 4 \max_{1\le i \le I, 1 \le j \le J} \|\br_{ij} - \br_{ij}^*\|  \backdefeq 4 M,$$ so that \begin{equation} \label{eq:b1} \abs{\|\br_u\|^2 - \|\br_u^*\|^2} \le 8 M.  \end{equation} The same bound holds for $\abs{\|\br_v\|^2 - \|\br^*_v\|^2}$. Next \begin{equation} \label{eq:b2}
\begin{split}
\abs{\innerprod{\br_u}{\br_v} - \innerprod{\br_u^*}{\br_v^*}} &= \abs{\innerprod{\br_u - \br_u^*}{\br_v} + \innerprod{\br_u^*}{\br_v - \br_v^*}}\\
& \le \|\br_u-\br_u^*\|\|\br_v\| + \|\br_u^*\|\|\br_v - \br_v^*\| \le 8M,
\end{split} 
\end{equation} as well.  Together \eqref{eq:b1} and \eqref{eq:b2} show that \begin{equation} \label{eq:gdiff}\|g-g^*\|_\infty \le C\cdot M. \end{equation} This implies that for a given curve in $\gamma$ in $(u,v)$-space, the length---measured in their respective metrics---of the curves $(\gamma,\br(\gamma))$ and $(\gamma,\br^*(\gamma))$ differ by at most $C\cdot M$, where $C$ has perhaps increased to account for the longest possible length of $\gamma$. Thus, taking $\gamma$ such that $(\gamma,\br(\gamma))$ is a geodesic connecting the origin and $(i,j,\br_{ij})$, we see $$D^*_{ij} \le \text{length}_{g^*}((\gamma,\br^*(\gamma))) \le \text{length}_g((\gamma,\br(\gamma))) + C M = D_{ij} + CM.$$ Taking instead $\gamma^*$ such that $(\gamma^*,\br^*(\gamma^*))$ is a geodesic connecting the origin to $(i,j,\br^*_{ij})$, we achieve the same bound with $D_{ij}$ and $D_{ij}^*$ reversed. Thus $$\abs{D_{ij} - D_{ij}^*} \le CM.$$ Taking the maximum on the left hand side completes the proof. \hfill $\square$\\

Next, corresponding to step (III) on the previous page, we prove an estimate relating two discrete immersions $\{\br_{ij},\bN_{ij}\}$ and $\{\br_{ij}^*,\bN^*_{ij}\}$ which are computed from the Lelieuvre equations with the same boundary data but different discrete curvature functions $K_{ij}$ and $K^*_{ij}$, respectively. 

\begin{lemma} \label{lem:2}
    Suppose that $K, K^* \le -1$ are two curvature functions, that give rise (via the initialization procedure described above, and the discrete Leliuevre formulas) to discrete surfaces $\{\br_{ij}, \bN_{ij}\}$ and $\{\br_{ij}^*, \bN_{ij}^*\}$ respectively. Further assume that there is $\mu > 0$ such that each of $$\zeta_{ij} =\|\bnu_{i,j-1} + \bnu_{i-1,j}\|, \,\,\,\,\,\ \zeta^*_{ij} =\|\bnu^*_{i,j-1} + \bnu^*_{i-1,j}\|$$ satisfy $\zeta_{ij}, \zeta^*_{ij} \ge \mu$ for all $i,j$. Then there is a constant $C>0$ independent of $K, K^*$ such that for all $i,j$, $$\|\br_{ij} - \br^*_{ij}\|,\,\,\|\bN_{ij} - \bN^*_{ij}\| \le C\|K-K\|_{\infty}.$$
\end{lemma}

\emph{Remark.} The condition that $\zeta_{ij}, \zeta^*_{ij}$ are bounded away from zero could likely be relaxed or rephrased, except possibly near singular edges where the normal vector flips orientation. Away from singular edges, it would require an extreme amount of curvature to accomplish $\bnu_{i+1,j} \approx -\bnu_{i,j+1}$, and thus force $\zeta_{ij} \approx 0$, so one may be able to derive a condition on $K$ which implies $\zeta_{ij} \ge \mu > 0$. In any of our figures, even in the presence of singular edges, $\zeta_{ij}$ remained fairly large. For example, for the sector in figure \ref{fig:branchPoint}, we have $\min_{ij}\zeta_{ij} \approx 1.6784$. Together with the fact that $\zeta_{ij} \le 2$, this indicates that $\bnu_{i+1,j}$ and $\bnu_{i,j+1}$ point in roughly the same direction, so that this assumption is not prohibitive.\\

\emph{Proof.} Recall, from the initialization procedure described in Section \ref{sec:amsler}, we choose the maximum distances $u_{\text{max}}$ and $v_{\text{max}}$ to extend the bounding lines for the Amsler-surface, and then fix $\|\br_{i+1,0}-\br_{i,0}\| \backdefeq \Delta u$, $\|\br_{0,j+1}-\br_{0,j}\| \backdefeq \Delta v$ where $I\Delta u = u_{\text{max}}$ and $J\Delta v = v_{\text{max}}.$ We proceed with this notation. 

Note that $\br_{ij}$ and $\br_{ij}^*$ are given the same values on the discrete boundaries $i=0$ and $j=0$. We first prove that the normals do not differ by too much along the boundaries. We prove this for the difference $\|\bN_{i,0} - \bN_{i,0}^*\|$; the same bound will hold \emph{mutatis mutandis} for $\|\bN_{0,j}-\bN_{0,j}^*\|$. Note that $\bN_{i,0}, \bN_{i,0}^*$ are given by \eqref{eq:normalBCs} with respective rotation angles $\delta_i, \delta_i^*$ given by \eqref{eq:rotationAngles} using either $K$ or $K^*$. In what follows, all norms are 2-norms unless otherwise specified, and $C$ is a constant which changes from line to line and depends on some of the ambient data, but \emph{not} upon $\|K-K^*\|_\infty$.   First note that because $\arcsin$ has derivative bounded by $1$, we have \begin{equation} \label{eq:deltaBound} \begin{split} \abs{\delta_i - \delta_i^*} &\le \Delta u\abs{ \abs{K_{i+1,0}K_{i,0}}^{-1/4} -  \abs{K^*_{i+1,0}K^*_{i,0}}^{-1/4}} \\ 
&\le C \cdot \Delta u \|K-K^*\|_{\infty}. 
\end{split} \end{equation} where $C$ is a constant depending on the maximum values of $K$ and $K^*$. Because the entries of a rotation matrix are smooth functions of the rotation angle, this gives $$\|\text{Rot}(\delta_i, (1,0,0)) - \text{Rot}(\delta_i^*, (1,0,0)) \| \le C\cdot \Delta u \|K-K^*\|_{\infty}.$$ Thus, using $R_i =\text{Rot}(\delta_i, (1,0,0))$ and $R^*_i = \text{Rot}(\delta^*_i, (1,0,0))$, we see 
\begin{align*}
\|\bN_{i,0} - \bN_{i,0}^*\| &= \left\| \left(\prod^{i-1}_{\ell = 0} R_\ell \right)\bN_{0,0} - \left( \prod^{i-1}_{\ell=0} R_\ell^*\right) \bN_{0,0} \right\|\\
&\le \left\| \prod^{i-1}_{\ell = 0} R_\ell- \prod^{i-1}_{\ell=0} R_\ell^* \right\|\\
&= \left\| \left(\prod^{i-1}_{\ell = 1} R_\ell \right)(R_0 - R^*_0) - \left( \prod^{i-1}_{\ell=1} R_\ell^* - \prod^{i-1}_{\ell=1} R_\ell\right)R^*_0 \right\| \\
&\le \|R_0 - R^*_0\| +  \left\| \prod^{i-1}_{\ell = 1} R_\ell- \prod^{i-1}_{\ell=0} R_\ell^* \right\|,
\end{align*} from which we arrive at \begin{equation} \label{eq:Ni0bound}
\|\bN_{i,0} - \bN^*_{i,0} \| \le \sum^{i-1}_{\ell=0} \|R_\ell - R^*_\ell\| \le C\cdot I \Delta u \|K-K^*\|_{\infty} = C \|K-K^*\|_\infty.
\end{equation} Similar steps lead to  \begin{equation} \label{eq:N0jbound}
\|\bN_{0,j} - \bN^*_{0,j}\| \le C \|K-K^*\|_\infty.
\end{equation} Next we need similar bounds for for any pair of $i = 1,\ldots, I$, $j = 1,\ldots, J$. From \eqref{eq:nuPreUpdate}, we have $$
\bN_{ij} = \abs{K_{ij}}^{1/4}\Big(\beta_{ij} (\bnu_{i,j-1} + \bnu_{i-1,j}) - \bnu_{i-1,j-1}\Big)
$$ where, defining $\zeta_{ij} =\|\bnu_{i,j-1} + \bnu_{i-1,j}\|$ and $\psi_{ij} = \innerprod{\bnu_{i,j-1} + \bnu_{i-1,j}}{\bnu_{i-1,j-1}} $, we have \begin{equation} \label{eq:cij}\beta_{ij} = \frac{\psi_{ij}+ \sqrt{\psi_{ij}^2 +\zeta_{ij}^2 (\abs{K_{ij}}^{-1/2} - \abs{K_{i-1,j-1}}^{-1/2})}}{\zeta_{ij}^2} =: F(\zeta_{ij},\psi_{ij}, K_{ij}, K_{i-1,j-1}) . \end{equation} Thus we have the bound \begin{equation} \label{eq:NijBound}\| \bN_{ij} - \bN^*_{ij}\| \le \text{I} + \text{II} \end{equation} where \begin{equation} \label{eq:IandII}\begin{split}\text{I} &= \|\beta_{ij} (\bnu_{i,j-1} + \bnu_{i-1,j}) - \bnu_{i-1,j-1}\|\abs{\abs{K_{ij}}^{1/4}-\abs{K^*_{ij}}^{1/4}}, \\ \text{II} &=\abs{K^*_{ij}}^{1/4}\|\beta_{ij} (\bnu_{i,j-1} + \bnu_{i-1,j}) - \beta_{ij}^*(\bnu^*_{i,j-1} + \bnu^*_{i-1,j}) - (\bnu_{i-1,j-1}-\bnu^*_{i-1,j-1}) \|. 
\end{split} \end{equation}  Since, we have assumed that $0 < \mu \le \zeta_{ij}, \zeta^*_{ij}$ and since $\nu,\nu^*$ are bounded, we see that $\beta_{ij}, \beta^*_{ij}$ remain bounded. In particular, this gives \begin{equation} \label{eq:Ibound} \text{I} \le C \|K-K^*\|_{\infty}.\end{equation} Next, we see \begin{equation}
\label{eq:IIbound}
\begin{split}
\text{II} &\le C \bigg(\abs{\beta_{ij}} \Big( \|\bnu_{i,j-1} - \bnu^*_{i,j-1}\| + \|\bnu_{i-1,j} - \bnu^*_{i-1,j}\|\Big) + \\ 
&\hspace{2cm} \|\bnu^*_{i,j-1} + \bnu^*_{i-1,j}\| \abs{\beta_{ij} - \beta^*_{ij}} + \|\bnu_{i,j-1} - \bnu^*_{i,j-1}\|\Big)
\end{split}
\end{equation} But the preceding observations---along with $\abs{K},\abs{K^*}\ge1$---show that the function $F$ defined in \eqref{eq:cij} is Lipschitz continuous. Thus \begin{align*} \abs{c_{ij} - c^*_{ij}} &= \abs{F(\zeta_{ij},\psi_{ij}, K_{ij}, K_{i-1,j-1}) - F(\zeta^*_{ij},\psi^*_{ij}, K^*_{ij}, K^*_{i-1,j-1})} \\
&\le C\|(\zeta_{ij},\psi_{ij}, K_{ij}, K_{i-1,j-1}) - (\zeta^*_{ij},\psi^*_{ij}, K^*_{ij}, K^*_{i-1,j-1})\|\\
&\le C(\abs{\zeta_{ij} - \zeta^*_{ij}} + \abs{\psi_{ij} - \psi_{ij}^*} + 2\|K-K^*\|_\infty). 
\end{align*} We notice that for any $i,j$, \begin{equation} \label{eq:nuBound} \begin{split} \|\bnu_{ij} - \bnu^*_{ij}\| &\le \abs{\abs{K_{ij}}^{-1/4} - \abs{K^*_{ij}}^{-1/4}}+ \abs{K^*_{ij}}^{-1/4} \|\bN_{ij} - \bN^*_{ij} \| \\
&\le C\|K-K^*\|_{\infty} + \|\bN_{ij} - \bN_{ij}^*\|. 
 \end{split} \end{equation} 
In particular, 
\begin{align*} \abs{\zeta_{ij}  - \zeta^*_{ij}} &\le \|\bnu_{i,j-1} - \bnu^*_{i,j-1}\| +  \|\bnu_{i-1,j} - \bnu^*_{i-1,j}\| \\ 
&\le C\|K-K^*\|_{\infty} + \|\bN_{i,j-1} - \bN_{i,j-1}^*\| + \|\bN_{i-1,j} - \bN_{i-1,j}^*\|
 \end{align*} and \begin{align*} 
\abs{\psi_{ij} - \psi_{ij}^*} &= \abs{\innerprod{\bnu_{i,j-1} + \bnu_{i-1,j}}{\bnu_{i-1,j-1}} -\innerprod{\bnu^*_{i,j-1} + \bnu^*_{i-1,j}}{\bnu_{i-1,j-1}^*}}\\
&\le \|\bnu_{i,j-1} - \bnu^*_{i,j-1} \| + \|\bnu_{i-1,j} - \bnu^*_{i-1,j} \| + \|\bnu_{i-1,j-1} - \bnu^*_{i-1,j-1} \|\\
&\le C\|K-K^*\|_{\infty} + \|\bN_{i,j-1} - \bN_{i,j-1}^*\| + \|\bN_{i-1,j} - \bN_{i-1,j}^*\| \\ &\hspace{8cm}+ \|\bN_{i-1,j-1} - \bN^*_{i-1,j-1}\|.
\end{align*} Inserting all these bounds into \eqref{eq:IIbound}, we see \begin{equation} 
\label{eq:IIbound2}
\begin{split}
\text{II} &\le C\Big( \|K-K^*\|_{\infty} + \|\bN_{i,j-1} - \bN_{i,j-1}^*\| \\ &\hspace{1.5cm}+ \|\bN_{i-1,j} - \bN_{i-1,j}^*\| + \|\bN_{i-1,j-1} - \bN^*_{i-1,j-1}\| \Big), \end{split}
\end{equation} and inserting \eqref{eq:Ibound} and \eqref{eq:IIbound2} into \eqref{eq:NijBound} results in \begin{equation} \label{eq:NijBound2} \begin{split}
\|\bN_{ij} - \bN^{*}_{ij}\| &\le C\Big( \|K-K^*\|_{\infty} + \|\bN_{i,j-1} - \bN_{i,j-1}^*\| \\ &\hspace{1.5cm}+ \|\bN_{i-1,j} - \bN_{i-1,j}^*\| + \|\bN_{i-1,j-1} - \bN^*_{i-1,j-1}\| \Big)
\end{split}
\end{equation} Recursively applying this until we hit the boundary and then using \eqref{eq:Ni0bound} and \eqref{eq:N0jbound} yields  \begin{equation} 
\label{eq:NijBoundFinal}
\|\bN_{ij} - \bN_{ij}^*\| \le C\left(\sum^{\max\{I,J\}}_{\ell=1} C^\ell\right)\|K-K^*\|_{\infty} = C\|K-K^*\|_{\infty}.
\end{equation} Finally, the discrete Lelieuvre formulas \eqref{eq:discreteLelieuvre} show that \begin{align*}
\| \br_{ij} - \br^*_{ij}\| &\le \|\br_{i-1,j} - \br^*_{i-1,j}\| +  \| \bnu_{ij}\times \bnu_{i-1,j} -  \bnu^*_{ij}\times \bnu^*_{i-1,j} \| \\
&\le \|\br_{i-1,j} - \br^*_{i-1,j}\| + \Big( \|\bnu_{ij} - \bnu^*_{ij}\| + \|\bnu_{i-1,j} - \bnu^*_{i-1,j}\|\Big)
\end{align*} whereupon \eqref{eq:nuBound} and \eqref{eq:NijBoundFinal} lead to
$$
\|\br_{ij} - \br_{ij}^*\| \le  \|\br_{i-1,j} - \br_{i-1,j}^*\| + C \|K-K^*\|_{\infty}. 
$$ Again, recursively applying the inequality until we reach the boundary, and using $\br_{i,0} = \br^*_{i,0}$, we have \begin{equation} 
\label{eq:rijBound}
\|\br_{ij} - \br^*_{ij}\| \le C\|K-K^*\|_{\infty}
\end{equation} which completes the proof. \hfill $\square$\\

With these lemmas, we can state and prove our convergence result. 

\begin{theorem}
    Assume that $h:[0,\infty)\to[0,\infty)$ is Lipschitz continuous, $\eps > 0$ is sufficiently small, and any discrete surfaces generated by the algorithm satisfy the assumption of lemma \ref{lem:2}. Then the iteration in Algorithm \ref{alg:IterationAlg} converges. 
\end{theorem}

\emph{Proof.} As stated at the beginning of the appendix, the iteration in Algorithm \ref{alg:IterationAlg} can be broken down into three steps, labeled (I), (II) and (III) above. Representing a discrete immersion $\br_{ij}$ and normal $\bN_{ij}$ as members of $\mathbb R^{3IJ}$, we can define the operator $T:\R^{3IJ} \times \R^{3IJ} \to \R^{3IJ} \times \R^{3IJ}$ so that the iteration can be written succinctly as $$\{\br^{(n+1)}_{ij}, \bN^{(n+1)}_{ij}\} = T\Big(\{\br^{(n)}_{ij}, \bN^{(n)}_{ij}\}\Big), \,\,\,\,\,\, n =0,1,2,\ldots.$$ We need to prove that this operator is contractive. We note that the underlying curvature $K_{ij}$ of the input surface is fixed in the definition of the operator $T$.

To do so, consider two discrete immersions and normal fields $\{\br_{ij}, \bN_{ij}\}, \{\br^{*}_{ij}, \bN^{*}_{ij}\}$. By lemma \ref{lem:1}, after step I, we arrive at discrete distance functions $D_{ij}$ and $D^*_{ij}$ such that $$\max_{ij}\abs{D_{ij} - D^*_{ij}} \le C\max_{ij}\|\br_{ij} - \br^*_{ij}\|$$ Then step II, we use these distance functions to update the curvature function; call the new values $K^{\text{new}}_{ij}, K^{*,\text{new}}_{ij}$. Then \begin{align*}\abs{K^{\text{new}}_{ij} - K^{*,\text{new}}_{ij}} &= \abs{(K_{ij}-\eps h(D_{ij})) - (K_{ij} - \eps h(D^*_{ij}))} \\ 
&= \eps \abs{h(D_{ij})-h(D^*_{ij})} \\ 
&\le \eps L \abs{D_{ij}-D^*_{ij}}\\
&\le \eps C \max_{ij} \|\br_{ij} - \br^*_{ij}\|,
\end{align*} where $L$ is the Lipschitz constant of $h$ (and is absorbed into $C$ on the ensuing line). Finally, we use these new discrete curvature functions in step III to generate new discrete immersions and normal fields (these are, by definition, the images of the old ones under $T$). By lemma \ref{lem:2}, we thus have \begin{align*}
   \max_{ij}\left \|T\Big(\{\br_{ij}, \bN_{ij}\}\Big) - T\Big(\{\br^{*}_{ij}, \bN^*_{ij}\}\Big)\right\| &\le C\|K^{\text{new}}_{ij}-K^{*,\text{new}}_{ij}\| \\
   &\le \eps C \max_{ij}\|\br_{ij} - \br^*_{ij}\|,
\end{align*} for some absolute constant $C > 0$. Taking $\eps < 1/C$, the operator $T$ will be a contraction, and thus the iteration will converge as a consequence of the Banach Fixed Point Theorem. \hfill $\square$ \\

\bibliographystyle{plain}
\bibliography{bibliography}

\end{document}

%% file: frontmatter.tex
\usepackage[top=1 in, bottom=1 in, left=1 in, right = 1 in]{geometry}
\usepackage{amsmath}
\usepackage{amsthm}
\usepackage{amssymb}
\usepackage{amsbsy}
\usepackage{setspace}
\usepackage{graphicx}
\usepackage{xcolor}
\definecolor{Hazel}{rgb}{0.3804, 0.4863, 0.3451}
\usepackage[colorlinks=true,citecolor=Hazel]{hyperref}
\usepackage[font=scriptsize,labelfont=bf]{caption}
\usepackage[font=scriptsize,labelfont=bf]{subcaption}
\graphicspath{ {images/} }
\usepackage{algorithm}
\usepackage{algpseudocode}
\usepackage[export]{adjustbox}  
\usepackage{tikz}

\newtheorem{theorem}{Theorem}[section]
\newtheorem{lemma}{Lemma}[section]
\newcounter{count}
\stepcounter{count}

\newcommand{\abs}[1]{\left| #1  \right|}

\newcommand{\B}[1]{\mathbf{#1}}
\newcommand{\R}{\mathbb{R}}
\newcommand{\Z}{\mathbb{Z}}

\newcommand{\br}{\mathbf{r}}
\newcommand{\bN}{\mathbf{N}}
\newcommand{\bv}{\mathbf{v}}
\newcommand{\bnu}{\boldsymbol{\nu}}
\newcommand{\innerprod}[2]{\langle #1, #2 \rangle}

\newcommand{\eps}{\varepsilon}

\newcommand*{\defeq}{\mathrel{\vcenter{\baselineskip0.5ex \lineskiplimit0pt
                     \hbox{\scriptsize.}\hbox{\scriptsize.}}}%
     		     =}
\newcommand*{\backdefeq}{=\mathrel{\vcenter{\baselineskip0.5ex \lineskiplimit0pt
                     \hbox{\scriptsize.}\hbox{\scriptsize.}}}%
     		     }

\def\XXint#1#2#3{{\setbox0=\hbox{$#1{#2#3}{\int}$ }
\vcenter{\hbox{$#2#3$ }}\kern-.6\wd0}}

\ifpdf
  \DeclareGraphicsExtensions{.eps,.pdf,.png,.jpg}
\else
  \DeclareGraphicsExtensions{.eps}
\fi

\usepackage{enumitem}
\setlist[enumerate]{leftmargin=.5in}
\setlist[itemize]{leftmargin=.5in}

\title{Discrete Differential Geometry for $C^{1,1}$ Hyperbolic Surfaces of Non-Constant Curvature}

\author{Christian Parkinson\thanks{Department of Mathematics and Department of Computational Mathematics, Science and Engineering, Michigan State University, East Lansing, MI 
  (\email{chparkin@msu.edu}).}
\and Shankar C. Venkataramani\thanks{Department of Mathematics, University of Arizona, Tucson, AZ 
  (\email{shankar@math.arizona.edu})}
}
\usepackage{amsopn}